\documentclass[hidelinks,onefignum,onetabnum]{siamart220329}
\usepackage{graphicx}
\usepackage{color}
\usepackage{amsmath,amstext,amssymb}
\usepackage{mathrsfs}
\usepackage{multirow}
\usepackage{algorithmic}
\usepackage{amsopn}
\usepackage{threeparttable}
\usepackage{hyperref}
\usepackage{subfigure}
{\bfseries}{\normalfont}
\newtheorem{assumption}{Assumption}[section]{\bfseries}{\normalfont}
\newtheorem{remark}{Remark}[section]{\bfseries}{\normalfont}
\headers{DFAP Algorithms for Nonconvex-Concave Minimax Problems}{Zi Xu,ZiQi Wang, JingJing Shen, and YuHong Dai}
\title{Derivative-free Alternating Projection Algorithms for General Nonconvex-Concave Minimax Problems\thanks{This work is supported by National Natural Science Foundation of China under the grant
			12071279.}}
\author{Zi Xu\thanks{Department of Mathematics, Shanghai University, Shanghai 200444, People’s Republic of China; Newtouch Center for Mathematics of Shanghai University, Shanghai 200444, People’s Republic of China
		(\email{xuzi@shu.edu.cn},\email{colwzq@shu.edu.cn},\email{sjj97@shu.edu.cn}).}
	\and ZiQi Wang\footnotemark[2]
	\and JingJing Shen\footnotemark[2]
	\and YuHong Dai\thanks{LSEC, ICMSEC, Academy of Mathematics and Systems Science,
		Chinese Academy of Sciences, Beijing 100190, China
		(\email{dyh@lsec.cc.ac.cn}). Corresponding author.}}
\ifpdf
\hypersetup{
	pdftitle={ZO-ARGP Algorithms for Nonconvex-Concave Minimax Problems},
	pdfauthor={Zi Xu, ZiQi Wang, JingJing Shen and YuHong Dai}
}
\fi

\begin{document}
	\maketitle
\begin{abstract}
In this paper, we study zeroth-order algorithms for nonconvex-concave minimax problems, which have attracted widely attention in machine learning, signal processing and many other fields in recent years. We propose a zeroth-order alternating randomized gradient projection (ZO-AGP) algorithm for smooth nonconvex-concave minimax problems, and its iteration complexity to obtain an $\varepsilon$-stationary point is bounded by $\mathcal{O}(\varepsilon^{-4})$, and the number of function value estimation is bounded by $\mathcal{O}(d_{x}+d_{y})$ per iteration. Moreover, we propose a zeroth-order block alternating randomized proximal gradient algorithm (ZO-BAPG) for solving block-wise nonsmooth nonconvex-concave minimax optimization problems, and the iteration complexity to obtain an $\varepsilon$-stationary point is bounded by $\mathcal{O}(\varepsilon^{-4})$ and the number of function value estimation per iteration is bounded by $\mathcal{O}(K d_{x}+d_{y})$. To the best of our knowledge, this is the first time that zeroth-order algorithms with iteration complexity gurantee are developed for solving both general smooth and block-wise nonsmooth nonconvex-concave minimax problems. Numerical results on data poisoning attack problem  and distributed nonconvex sparse principal component analysis problem validate the efficiency of the proposed algorithms.
\end{abstract}
	
	\begin{keywords}
		nonconvex-concave minimax problem, zeroth-order algorithm, alternating randomized gradient projection algorithm, alternating randomized proximal gradient algorithm, complexity analysis, machine learning
	\end{keywords}
	
	\begin{MSCcodes}
		90C47, 90C26, 90C30
	\end{MSCcodes}
	
	\section{Introduction}
	We consider the following deterministic minimax optimization problem:
	\begin{equation}\label{problem}
		\min_{x \in \mathcal{X}} \max_{y \in \mathcal{Y}} f\left(x,y\right),\tag{P}
	\end{equation}
	where $f: \mathcal{X} \times \mathcal{Y} \rightarrow \mathbb{R}$ is a continuously differentiable function, $\mathcal{X} \subset \mathbb{R}^{d_{x}}$ and $\mathcal{Y} \subset \mathbb{R}^{d_{y}}$ are nonempty compact convex sets, $d_{x}$ and $d_{y}$ are dimension of $\mathcal{X}$ and $\mathcal{Y}$ respectively. In this paper, we study nonconvex-concave minimax problems, i.e., $f(x, y)$ is nonconvex w.r.t. $x$ and concave w.r.t. $y$. Moreover, we focus on the zeroth-order optimization algorithm to solve general nonconvex-concave minimax problems with a black-box setting, where the desired optimizer cannot access the gradients of the objective function but may query its values. Recently, this problem has attracted increasing attention in solving machine learning or deep learning problems, and many other research fields, e.g., the adversarial attack to deep neural networks(DNNs), reinforcement  learning, robust training, hyperparameter tuning~\cite{chen2017zoo,finlay2021scaleable,snoek2012practical}.
	
	There are some existing works that focus on zeroth order algorithms for solving general minimax optimization problems. For example,  derivative-free algorithms, Bayesian optimization algorithms and evolutionary algorithms are proposed for minimax optimization problems~\cite{ADSHO19, BN10,MW20,PBH19}. However, the works mentioned above do not provide any oracle complexity analysis.  For convex-concave minimax optimization problems, there are some exisiting algorithms. For example, Roy et al. \cite{roy2019online} study zeroth-order Frank-Wolfe algorithms for strongly convex-strongly concave minimax optimization problems and provide non-asymptotic oracle complexity analysis. Beznosikov et al. \cite{beznosikov2020gradient} present a zeroth-order saddle-point algorithm (zoSPA) with the total complexity  of $\mathcal{O}\left( \varepsilon ^{-2} \right)$. Sadiev et al. \cite{sadiev2020zeroth} propose zeroth-order analogues of Mirror-Descent and Mirror-Prox methods for stochastic convex-concave and convex-strongly concave minimax problems, with the total complexity  of $\mathcal{O}\left( \varepsilon ^{-2} \right)$ and $\mathcal{O}\left( \varepsilon ^{-1} \right)$ respectively. Maheshwari et al. \cite{maheshwari2022zeroth} propose an optimistic gradient descent ascent with random reshuffling (OGDA-RR) algorithm for solving convex-concave minimax problem with the total complexity  of $\mathcal{O}\left( \varepsilon ^{-1} \right)$.
	
	A few zeroth-order algorithms have been proposed for nonconvex-strongly concave minimax optimization problems.  Liu et al. \cite{liu2020min} propose an alternating projected stochastic gradient descent-ascent method, integrating a zeroth-order (ZO) gradient estimator, referred to as ZO-Min-Max algorithm and it can find an $\varepsilon$-stationary point with the total complexity of  $\mathcal{O}\left(\left(d_{x}+d_{y}\right) \varepsilon^{-4}\right)$.  Wang et al.~\cite{wang2020zerothorder} propose a zeroth-order gradient descent ascent (ZO-GDA) algorithm and a zeroth-order gradient descent multi-step ascent (ZO-GDMSA) algorithm for deterministic nonconvex-strongly concave minimax optimization problems, and the total number of calls to the (deterministic) zeroth-order oracle is $\mathcal{O}\left(\kappa^{5}\left(d_{x}+d_{y}\right) \epsilon^{-2}\right)$ and $\mathcal{O}\left(\kappa\left(d_{x}+\kappa d_{y} \log \left(\epsilon^{-1}\right)\right) \epsilon^{-2}\right)$ respectively. The stochastic version of the two algorithms are also proposed, i.e., ZO-SGDA and ZO-SGDMSA, and the corresponding total number of calls to the stochastic zeroth-order oracle is $\mathcal{O}\left(\kappa^{5}\left(d_{x}+d_{y}\right) \epsilon^{-4}\right)$ and $\mathcal{O}\left(\kappa\left(d_{x}+\kappa d_{y}\log \left(\epsilon^{-1}\right)\right) \epsilon^{-4}\right)$ respectively. Xu et al. \cite{xu2020enhanced} propose a zeroth-order variance reduced gradient descent ascent (ZO-VRGDA) algorithm, which achieves the total complexity of $\mathcal{O}\left( \kappa^{3}\left( d_{x}+d_{y}\right) \varepsilon^{-3}\right)$ under stochastic setting. Huang et al.~\cite{huang2020accelerated} propose an accelerated zeroth-order momentum descent ascent (Acc-ZOMDA) method and it can find an $\varepsilon$-stationary point with the total complexity of $\mathcal{O}\left(\kappa^{3}\left(d_{x}+d_{y}\right)^{3 / 2} \epsilon^{-3}\right)$.
	
	Shen et al. \cite{shen2022zeroth} propose a new zeroth-order alternating randomized gradient projection algorithm to solve smooth nonconvex-linear minimax problems and its iteration complexity to find an $\epsilon$-first-order Nash equilibrium is $\mathcal{O}(\varepsilon^{-3})$ and the number of function value estimation per iteration is bounded by $\mathcal{O}(\varepsilon^{-2})$. 
	
    There are also some related works to solve nonconvex-nonconcave minimax problems with special structures. A common assumption is that $f(x,y)$ satisfies the Polyak-\L ojasiewicz (PL) condition with respect to $y$. Xu et al. \cite{xu2023zeroth} propose a zeroth-order alternating gradient descent ascent (ZO-AGDA) algorithm and a zeroth-order variance reduced alternating gradient descent ascent (ZO-VRAGDA) algorithm for solving nonconvex-PL minimax problem under the deterministic and the stochastic setting, respectively. The corresponding comlexity number of of ZO-AGDA and ZO-VRAGDA algorithm to obtain an $\epsilon$-stationary point is upper bounded by $\mathcal{O}(\varepsilon^{-2})$ and $\mathcal{O}(\varepsilon^{-3})$, respectively.
	
	To the best of our knowledge, instead of nonconvex-strongly concave minimax optimization problems, there is no known existing algorithm with the complexity gurantee for solving general zeroth-order nonconvex-concave minimax problems.
	
	\subsection{Related Works}
	In this subsection, we give a brief review on the first order methods for solving minimax optimization problems. For solving convex-concave minimax optimization problems, there are many existing works in the early literature. Since we focus on nonconvex minimax problems, we do not attempt to survey it in this paper and refer to \cite{nemirovski2004prox,nesterov2007dual,mokhtari2020unified,lin2020near,chen2014optimal,chen2017accelerated,gidel2018variational,juditsky2016solving,lan2016iteration,mertikopoulos2018mirror, monteiro2011complexity, ouyang2019lower} for related results.
	
	For nonconvex-strongly concave minimax optimization problems, various algorithms have been proposed in \cite{jin2019minmax,lin2020gradient,lu2020hybrid,rafique2018non}, and all of them can achieve the gradient complexity of $\tilde{\mathcal{O}}\left( \kappa_y^2\varepsilon ^{-2} \right)$ in terms of stationary point of $\Phi (\cdot) = \max_{y\in \mathcal{Y}} f(\cdot, y)$ (when $\mathcal{X}=\mathbb{R}^n$), or that of $f$,  where $\kappa_y$ is the condition number for $f(x,\cdot)$. Lin et al. \cite{lin2020near} propose an accelerated algorithm which can improve the gradient complexity bound to $\tilde{\mathcal{O}}\left( \sqrt{\kappa_y}\varepsilon^{-2} \right)$. 
	
	For general nonconvex-concave (but not strongly concave) minimax problem, there are two types of algorithms, i.e., nested-loop algorithms and single-loop algorithms. One intensively studied type is nested-loop algorithms. Various algorithms of this type have been proposed in \cite{rafique2018non,nouiehed2019solving,thekumparampil2019efficient,kong2019accelerated,ostrovskii2020efficient,yang2020catalyst}. Lin et al. \cite{lin2020near} propose a class of accelerated algorithms for smooth nonconvex-concave minimax problems with the complexity bound of $\tilde{\mathcal{O}}\left( \varepsilon ^{-2.5} \right)$ which owns the best iteration complexity till now. 
	On the other hand, fewer studies have been directed to single-loop algorithms for nonconvex-concave minimax problems. One such method is the gradient descent-ascent (GDA) method, which performs a gradient descent step on $x$ and a gradient ascent step on $y$ simultaneously at each iteration. Many improved GDA algorithms have been proposed and analyzed in \cite{Letcher,Chambolle,Daskalakis17,Daskalakis18,gidel2018variational,Ho2016,lin2020gradient}. Jin et al. \cite{jin2019minmax} propose a GDmax algorithm with iteration complexity of $\tilde{\mathcal{O}}\left( \varepsilon ^{-6} \right) $. Moreover, Lu et al.~\cite{lu2020hybrid} propose the hybrid block successive approximation (HiBSA) algorithm, which can obtain an $\varepsilon$-stationary point of $f(x,y)$ in $\tilde{\mathcal{O}}\left( \varepsilon ^{-4} \right)$ iterations when $\mathcal{X}$ and $\mathcal{Y}$ are  convex compact sets. 
	Xu et al. \cite{xu2020unified} propose a unified single-loop alternating gradient projection (AGP) algorithm for solving nonconvex-(strongly) concave and (strongly) convex-nonconcave minimax problems, which can find an $\varepsilon$-stationary point with the gradient complexity of $\mathcal{O}\left( \varepsilon ^{-4} \right)$. Zhang et al. \cite{zhang2020single} propose a smoothed GDA algorithm which achieves the same iteration complexity for general nonconvex-concave minimax problems. Moreover, different stopping rules have been adopted in these existing works, e.g., \cite{lin2020gradient,lu2020hybrid,nouiehed2019solving}.

	Moreover, some works focus on the stochastic nonconvex-strongly concave minimax problems. For example,  Luo et al. \cite{luo2020stochastic}  propose a stochastic recursive gradient descent ascent (SREDA) algorithm, and Xu et al. \cite{xu2020enhanced} propose an enhanced SREDA algorithm, and the gradient complexity of both of them is $\tilde{\mathcal{O}}\left( \kappa^{3}\varepsilon^{-3} \right)$. Furthermore, several works study the lower complexity bounds to find approximate stationary points for nonconvex-strongly concave minimax problems, please refer to \cite{zhang2021complexity} and the references therein. 
	We refer to two recent review papers \cite{LMW19,liu2020primer} for literature on zeroth-order optimization methods for general argmin-type optimization problems, which we do not attempt to review in this paper.
	
	\subsection{Contributions}
	We propose a zeroth-order alternating randomized gradient projection (ZO-AGP) algorithm for smooth nonconvex-concave minimax problems. To obtain an $\varepsilon$-stationary point,  the iteration complexity of ZO-AGP algorithm  is proved to be bounded by $\mathcal{O}(\varepsilon^{-4})$ with the number of function value estimation per iteration being bounded by $\mathcal{O}(d_{x}+d_{y})$. To the best of our knowledge, this is the first time that a zeroth-order algorithm with iteration complexity gurantee is developed for general smooth nonconvex-concave minimax problems.
	
	
	Furthermore, we propose a zeroth-order block alternating randomized proximal gradient algorithm (ZO-BAPG) for solving block-wise nonsmooth nonconvex-concave minimax optimization problems, which owns the same iteration complexity with that of  ZO-AGP algorithm, i.e,  $\mathcal{O}(\varepsilon^{-4})$,  and the number of function value estimation per iteration is bounded by $\mathcal{O}(K d_{x}+d_{y})$. To the best of our knowledge, this is also the first time that a zeroth-order algorithm with iteration complexity gurantee is developed for solving general block-wise nonsmooth nonconvex-concave minimax problems.

	{\bfseries Notations}. We use $\mathbb{R}^{d_{x}}$ to denote the space of $d_{x}$ dimension real valued vectors. $\left\langle x,y \right\rangle$ denotes the inner product of two vectors of $x$ and $y$. $\left\| \cdot \right\|$ denotes the Euclidean norm. $\operatorname{Proj}_{\mathcal{X}}$ denotes the projection operator to the set $\mathcal{X}$. $\left[q\right]$ denotes the set $\{1,\cdots,q \}$, when $q$ is a positive integer. $\mathbb{E}_{u} \left(\cdot\right)$ denotes that we take expectation over random variable $u$. $\partial h$ denotes sub-gradient of function $h$.

	The remainder of this paper is organized as follows.  In Section \ref{section ZOAGP}, we propose a zeroth-order alternating gradient projection algorithm (ZO-AGP) for solving nonconvex-concave deterministic minimax optimization problems. In Section \ref{section Analysis_block}, we further propose a zeroth-order block alternating proximal gradient (ZO-BAPG) algorithm for solving block-wise nonconvex-concave deterministic minimax optimization problems, and  iteration complexity for both two proposed algorithms are also analyzed. Some numerical results on data poisoning attack problem and  distributed nonconvex sparse principal component analysis problem are shown in Section \ref{Experiments}. Some conclusions are made in the last section.

	\section{Zeroth-Order Algorithms for Nonconvex-Concave Minimax Problems}
	\subsection{A Zeroth-Order Alternating Gradient Projection Algorithm}\label{section ZOAGP}
	In this subsection, we propose a zeroth-order alternating  gradient projection (ZO-AGP) algorithm for solving \eqref{problem}.
	Before proposing our algorithm, 
	we first introduce a gradient estimator \cite{nemirovskij1983problem}, 
	which can estimate the gradient of function $f\left(x,y\right): \mathbb{R}^{d_{x}} \times \mathbb{R}^{d_{y}} \rightarrow \mathbb{R}$. The gradient estimator can be defined as
	\begin{align}
		\widehat{\nabla}_{x} f\left(x,y\right) &= \sum_{i=1}^{d_x} \frac{\left[f\left(x+\mu_{1} u_{i},y\right) - f(x,y)\right]}{\mu_{1}} u_{i}, \label{x_esti}\\
		\widehat{\nabla}_{y} f\left(x,y\right) &= \sum_{i=1}^{d_y} \frac{\left[f(x,y+\mu_{2} v_{i})-f(x,y)\right]}{\mu_{2}} v_{i}, \label{y_esti}
	\end{align}
	where  $\mu_{1},\ \mu_{2} > 0$ are smoothing parameters, $\{u_{i}\}_{i=1}^{d_x}$ is a standard basis in $\mathbb{R}^{d_{x}}$ and $\{v_{i}\}_{i=1}^{d_{y}}$ is a standard basis in $\mathbb{R}^{d_{y}}$.

	Next, we propose a zeroth-order alternating gradient projection (ZO-AGP) algorithm based on the idea of the AGP algorithm \cite{xu2020unified}. Note that each iteration of the AGP algorithm consists of two gradient projection steps for updating both $x$ and $y$. Instead of the original function $f(x, y)$, AGP uses the gradient of a regularized version of the original function, i.e., $ \tilde{f}_t(x, y) = f(x, y)  -\frac{\lambda_t}{2}\|y\|^2$
	where $\eta_t \ge 0$ are regularization parameters. Similarly, each iteration of the novel proposed ZO-AGP algorithm conducts two ``gradient" projection steps for updating both $x$ and $y$, where the ``gradient" is obtained by gradient estimators shown in \eqref{x_esti} and \eqref{y_esti}. More detailedly, it updates $x_t$ by minimizing a linearized approximation of $\tilde{f}_t \left( x,  y_t\right)$ where the gradient at point $\left( x_t,y_t \right)$ is replaced by a gradient estimators shown in  \eqref{x_esti}, i.e.,
	\begin{align}
		x_{t+1}&={\arg\min}_{x\in \mathcal{X}}\langle  \widehat{\nabla}_x \tilde{f}_t\left( x_t,y_t \right) ,x-x_t \rangle  +\tfrac{1}{2\alpha_t}\|x-x_t \| ^2\nonumber\\
		&= \operatorname{proj}_\mathcal{X} \left(x_t - \alpha_t  \widehat{\nabla}_x \tilde{f}_t ( x_t,y_t ) \right), \label{alg:x}
	\end{align}
	where $\operatorname{proj}_\mathcal{X}$ is the projection operator onto $\mathcal{X}$, $ \alpha_t  > 0$ denotes a stepsize parameter, and
	\begin{align}
		\widehat{\nabla}_{x}  \tilde{f}_t\left(x,y\right) = \sum_{i=1}^{d_x} \frac{\left[f\left(x+\mu_{1,t} u_{i},y\right) - f(x,y)\right]}{\mu_{1,t}} u_{i}, \label{x_esti_t}
	\end{align}
	where $\{\mu_{1,t}\}$ is a decreasing sequence.
	Similarly, it updates $y_t$ by maximizing a linearized approximation of $\tilde{f}\left( x_{t+1}, y\right)$ minus some regularized term where the gradient at point $\left( x_{t+1},y_t \right)$ is replaced by a gradient estimator shown in  \eqref{y_esti}, i.e.,
	\begin{align}
		y_{t+1}&={\arg\max}_{y\in \mathcal{Y}}\langle \widehat{\nabla}_y\tilde{f}_t\left( x_{t+1},y_t \right),y-y_t \rangle  -\tfrac{1}{2\beta}\| y-y_t \| ^2\nonumber\\
		&= \operatorname{proj}_{\mathcal{Y}}\left(y_{t}+\beta \left(\widehat{\nabla}_{y} f\left(x_{t+1}, y_{t}\right)-\lambda_ty_{t}\right)\right),\label{alg:y}
	\end{align}
	where  $\operatorname{proj}_\mathcal{Y}$ is the projection operator onto $\mathcal{Y}$, $ \beta  > 0$ denotes a stepsize parameter, and
	\begin{equation}
		\widehat{\nabla}_{y} \tilde{f}_t\left(x,y\right)=\widehat{\nabla}_{y} f\left(x, y\right)-\lambda_ty =  \sum_{i=1}^{d_y} \frac{\left[f(x,y+\mu_{2} v_{i})-f(x,y)\right]}{\mu_{2}} v_{i}-\lambda_ty. \label{y_esti_t}
	\end{equation}
	The proposed ZO-AGP algorithm is formally stated in Algorithm \ref{algo1}.

\begin{remark}
	Note that the Algorithm \ref{algo1} differs from the ZO-Min-Max algorithm proposed in \cite{liu2020min} in two ways. On the one hand, the ZO-Min-Max algorithm uses uniformly distributed random vectors on the unit sphere to calculate the zero-order gradient estimator, while our algorithm uses a finite-difference stochastic gradient estimator which can yields to a better complexity bound. On the other hand, the ZO-Min-Max algorithm is used to solve nonconvex-strongly concave minimax problems, whereas Algorithm \ref{algo1} can be used to solve more general nonconvex-concave minimax problems.
\end{remark}

		\begin{algorithm}
			\caption{(ZO-AGP Algorithm)}
			\label{algo1}
			\begin{algorithmic}
				\STATE{\textbf{Step 1}:  Input $x_1,y_1,\beta,\alpha_1, \lambda_1$; Set $t=1$.}
				\STATE{\textbf{Step 2}:~Calculate $\alpha_t$, and perform the following update for $x_t$:  	\begin{equation}\label{algolstep_x}
						x_{t+1}=\operatorname{proj}_{\mathcal{X}}\left(x_{t}-\alpha_{t} \widehat{\nabla}_x \tilde{f}_t\left(x_{t}, y_{t}\right)\right).
				\end{equation}}
				\STATE{\textbf{Step 3}:~Calculate $\eta_t$, and perform the following update for $y_t$:
					\begin{equation}\label{algolstep_y}
						y_{t+1}=\operatorname{proj}_{\mathcal{Y}}\left(y_{t}+\beta \left(\widehat{\nabla}_{y} f\left(x_{t+1}, y_{t}\right)-\lambda_ty_{t}\right)\right).
				\end{equation}}
				\STATE{\textbf{Step 4}:~If converges, stop; otherwise, set $t=t+1, $ go to Step 2.}
			\end{algorithmic}
		\end{algorithm}
		
		To analyze the convergence of Algorithm \ref{algo1}, we define the stationarity gap as the termination criterion as follows.
		
		\begin{definition}
			At each iteration of Algorithm \ref{algo1}, the stationarity gap for problem \eqref{problem} is defined as:
			\begin{equation}
				\nabla\mathbf{G}\left(x_{t},y_{t}\right) =
				\begin{bmatrix}
					\frac{1}{\alpha_t}\left(x_{t} - \operatorname{Proj}_{\mathcal{X}} \left(x_{t}- \alpha_{t} \nabla_{x}f\left(x_{t},y_{t}\right)\right)\right)\\
					\frac{1}{\beta} \left(y_{t} - \operatorname{Proj}_{\mathcal{Y}} \left(y_{t}+\beta\nabla_{y}f\left(x_{t},y_{t}\right)\right)\right)
				\end{bmatrix}.
			\end{equation}
		\end{definition}
		
		\begin{definition}
			At each iteration of Algorithm \ref{algo1}, we can also define the stationarity gap as:
			\begin{equation}\label{conver_block}
				\nabla\tilde{\mathbf{G}}\left(x_{t},y_{t}\right) =
				\begin{bmatrix}
					\frac{1}{\alpha_{t}} \left(x_{t}- \operatorname{Proj}_{\mathcal{X}} \left(x_{t}- \alpha_t \nabla_{x}\tilde{f}_t\left(x_{t},y_{t}\right)\right)\right)\\
					\frac{1}{\beta} \left(y_{t} - \operatorname{Proj}_{\mathcal{Y}} \left(y_{t}+\beta\nabla_{y}\tilde{f}_t\left(x_{t},y_{t}\right)\right)\right)
				\end{bmatrix}.
			\end{equation}
			For simplicity, we denote $\nabla\tilde{\mathbf{G}}_{t}=\nabla\tilde{\mathbf{G}}\left(x_{t},y_{t}\right)$,
			\begin{align*}
				\left(\nabla\tilde{\mathbf{G}}_{t}\right)_{x}&=\frac{1}{\alpha_{t}} \left(x_{t} - \operatorname{Proj}_{\mathcal{X}} \left(x_{t}- \alpha_t \nabla_{x}\tilde{f}_t\left(x_{t},y_{t}\right)\right)\right),\\
				\left(\nabla\tilde{\mathbf{G}}_{t}\right)_{y}&=\frac{1}{\beta} \left(y_{t} - \operatorname{Proj}_{\mathcal{Y}} \left(y_{t}+\beta\nabla_{y}\tilde{f}_t\left(x_{t},y_{t}\right)\right)\right).
			\end{align*}
		\end{definition}
		

		\subsection{Complexity analysis}\label{section Analysis_b}
		In this subsection, we prove the iteration complexity of 
		ZO-AGP algorithm for solving smooth nonconvex-concave deterministic minimax optimization problems. 
		
		Firstly, we need to make the following assumption about the smoothness of $f(x,y)$.
		
		\begin{assumption}\label{a4}
			$f\left(x,y\right)$ has Lipschitz continuous gradients, i.e., there exist $L_{x}$, $L_{y}$ such that,
			\begin{align*}
				\left\|\nabla_{x}f\left(x_{1},y\right) - \nabla_{x}f\left(x_{2},y\right)\right\| & \leq L_{x} \left\|x_{1}-x_{2}\right\|, \quad \forall x_{1},x_{2}\in\mathcal{X}, \forall y\in \mathcal{Y}, \\
				\left\|\nabla_{y}f\left(x_{1},y\right) - \nabla_{y}f\left(x_{2},y\right)\right\| &\leq L_{y} \left\|x_{1}-x_{2}\right\|, \quad \forall x_{1},x_{2}\in\mathcal{X},\forall y\in \mathcal{Y}, \\
				\left\|\nabla_{y}f\left(x,y_{1}\right) - \nabla_{y}f\left(x,y_{2}\right)\right\| &\leq L_{y} \left\|y_{1}-y_{2}\right\|, \quad \forall x\in \mathcal{X}, \forall y_{1},y_{2}\in\mathcal{Y}.
			\end{align*}
		\end{assumption}
		
		
		We also make the following assumption about the parameters $\lambda_t$.
		\begin{assumption}\label{ass:2}
			$\{\lambda_t\}$ is a nonnegative monotonically decreasing sequence.
		\end{assumption}

		Next, we provide an upper bound on the variance gradient estimators as follows.
		
		\begin{lemma}
			\label{varboundlemma_b}
			Suppose that Assumption \ref{a4} holds.  Then we have  
			\begin{align}
				\left\|\widehat{\nabla}_x \tilde{f}_t \left(x,y\right) - \nabla_{x}f\left(x,y\right)\right\|^{2} &\leq \frac{d_xL_x^2\mu_{1,t}^2}{4},\\
				\left\|\widehat{\nabla}_{y}f\left(x,y\right) - \nabla_{y}f\left(x,y\right)\right\|^{2} &\leq \frac{d_yL_y^2\mu_2^2}{4}.
			\end{align}
		\end{lemma}
		
		\begin{proof}
			By Assumption \ref{a4},  for all $i = 1, \cdots, d_x$, we have	
			\begin{equation}
				|f(x+\mu_{1,t}u_i,y)-f(x,y)-\langle\nabla_{x}f(x,y),\mu_{1,t}u_i\rangle|\le\frac{L_x}{2}\|\mu_{1,t}u_i\|^2.\label{lem2.3:3}
			\end{equation}
		Since $\{u_{i}\}_{i=1}^{d_x}$ is a standard basis in $\mathbb{R}^{d_{x}}$, we get $\|u_i\|=1$.
		Then \eqref{lem2.3:3} implies that
			\begin{equation*}
				\left|\frac{f(x+\mu_{1,t}u_i,y)-f(x,y)}{\mu_{1,t}}-\langle\nabla_{x}f(x,y),u_i\rangle\right|\le\frac{L_x\mu_{1,t}}{2}.
			\end{equation*}
			By the definition of $\widehat{\nabla}_{x}\tilde{f}_t(x,y)$ in \eqref{x_esti_t}, we have
			\begin{align*}
				\|\widehat{\nabla}_x \tilde{f}_t(x,y)-\nabla_{x}f(x,y) \|^2&=\sum_{i=1}^{d_x}\left|\frac{f(x+\mu_{1,t}u_i,y)-f(x,y)}{\mu_{1,t}}-\langle\nabla_{x}f(x,y),u_i\rangle\right|^2\nonumber\\
				&\le\frac{d_xL_x^2\mu_{1,t}^2}{4}.
			\end{align*}
			Similarly, we can prove that
			\begin{align*}
				\|\widehat{\nabla}_{y}f(x,y)-\nabla_{y}f(x,y) \|^2
				\le\frac{d_yL_y^2\mu_2^2}{4}.
			\end{align*}
		\end{proof}
		
		Since that $f\left(x,y\right)$ is a continuous function, and $\mathcal{X}$, $\mathcal{Y}$ are  nonempty compact convex sets, $f(x,y)$ is bounded, i.e., there exist finite numbers $\bar{f}$ and $\underline{f}$, such that $\underline{f}\leq f(x,y) \leq \bar{f}$. 
		Now, we are ready to show the descent property w.r.t. $x$ as follows.
		\begin{lemma}\label{dec_x_blocklemma}
			Suppose that Assumption \ref{a4} holds. 
			Consider the sequence $\{\left(x_{t},y_t\right)\}$ generated by Algorithm \ref{algo1}. Then we have
			\begin{equation}\label{decentx_block}
				f(x_{t+1},y_{t})\leq f(x_{t},y_{t})- \left(\frac{1}{\alpha_{t}} - \frac{L_x}{2}-\frac{1}{2\beta}\right)\left\|x_{t+1} - x_{t}\right\|^{2} +\frac{\beta d_xL_x^2\mu_{1,t}^2}{8}.
			\end{equation}
		\end{lemma}
		
		\begin{proof}
			By Assumption \ref{a4}, $f(x,y)$ has Lipschtiz continuous gradient w.r.t. $x$, which implies that
			\begin{align}\label{l6(0)_block}
				f(x_{t+1},y_{t})
				&\leq f(x_{t},y_{t}) + \left\langle \nabla_{x} f(x_{t},y_{t}),x_{t+1}-x_{t} \right\rangle + \frac{L_{x}}{2}\|x_{t+1}-x_{t}\|^{2}\notag\\
				&= f(x_{t},y_{t}) + \left\langle\widehat{\nabla}_x \tilde{f}_t(x_{t},y_{t}),x_{t+1} - x_{t}\right\rangle + \frac{L_{x}}{2}\|x_{t+1}-x_{t}\|^{2}\notag\\
				&\quad +\left\langle\nabla_{x}f(x_{t},y_{t}) - \widehat{\nabla}_x \tilde{f}_t(x_{t},y_{t}),x_{t+1} - x_{t}\right\rangle.
			\end{align}
			Firstly, we estimate the second term in the r.h.s of \eqref{l6(0)_block}. From the optimality of $x_{t+1}$'s update in \eqref{algolstep_x} of  Algorithm \ref{algo1},  $\forall x\in\mathcal{X}$, we obtain
			\begin{equation}\label{l6_2nd1_block}
				\left\langle x_{t+1} - x_{t} +
				\alpha_{t}\widehat{\nabla}_x \tilde{f}_t(x_{t},y_{t}), x-x_{t+1}\right\rangle \geq 0.
			\end{equation}
			By setting $x=x_{t}$ and rearranging the terms, we have
			\begin{equation}\label{l6_2nd2_block}
				\left\langle\widehat{\nabla}_x \tilde{f}_t(x_{t},y_{t}),x_{t+1} - x_{t}\right\rangle \leq - \frac{1}{\alpha_{t}}\|x_{t+1}-x_{t}\|^{2}.
			\end{equation}
			On the other hand, we estimate the last term in the right side of \eqref{l6(0)_block}. By the Cauchy-Schwarz inequality and Lemma \ref{varboundlemma_b}, we have
			\begin{align}
				\left\langle\nabla_{x}f(x_{t},y_{t}) - \widehat{\nabla}_x \tilde{f}_t(x_{t},y_{t}), x_{t+1}-x_{t}\right\rangle
				\le\frac{1}{2\beta}\|x_{t+1}-x_{t}\|^{2}+\frac{\beta d_xL_x^2\mu_{1,t}^2}{8}.\label{l6_4th2_b}
			\end{align}
			The proof is then completed by plugging \eqref{l6_2nd2_block} and \eqref{l6_4th2_b} into \eqref{l6(0)_block}.
		\end{proof}

		Before we show the ascent quantity in $y$-maximization, we need the following lemma to show the recurrence of the distance between successive two iterations on $y$.

		\begin{lemma}\label{dislemma_block}
			Suppose that Assumption \ref{a4} holds. Assume $\beta \leq \frac{1}{L_{y}+2\lambda_{1}}$. Then, for the sequence $\{(x_{t},y_{t})\}$ generated by Algorithm \ref{algo1}, for any $t\geq 3$, we have
			\begin{align}\label{disy_block}
				&\frac{16}{\beta^{2}\lambda_{t+1}}\left\|y_{t+1}-y_{t}\right\|^{2}
				-\frac{16}{\beta}\left(\frac{\lambda_{t-1}}{\lambda_{t}}-1\right)\left\|y_{t+1}\right\|^{2}\notag\\
				\leq&\frac{16}{\beta^{2}\lambda_{t}}\left\|y_{t}-y_{t-1}\right\|^{2}
				-\frac{16}{\beta}\left(\frac{\lambda_{t-2}}{\lambda_{t-1}}-1\right)\left\|y_{t}\right\|^{2}
				+\frac{16}{\beta^{2}}\left(\frac{1}{\lambda_{t+1}}-\frac{1}{\lambda_{t}}\right) \left\|y_{t+1}-y_{t}\right\|^{2} \notag\\
				&
				-\frac{12}{\beta}\left\|y_{t}-y_{t-1}\right\|^{2}
				+\frac{2}{\beta}\left\|y_{t+1}-y_{t}\right\|^{2}
				+\frac{384}{\beta\lambda_{t}^{2}}L_{y}^{2}\left\|x_{t+1}-x_{t}\right\|^{2}
				\nonumber\\
				&+\frac{16}{\beta}\left(\frac{\lambda_{t-2}}{\lambda_{t-1}}-\frac{\lambda_{t-1}}{\lambda_{t}}\right)\left\|y_{t}\right\|^{2}+\left(\frac{32}{\lambda_{t}}+\frac{256}{\beta\lambda_{t}^{2}}\right)d_{y}L_{y}^{2}\mu_{2}^{2}.
			\end{align}
		\end{lemma}
		\begin{proof}
			Firstly, from the optimality condition of \eqref{algolstep_y} at $t+1$th and $t$th iteration, respectively, we have
			\begin{align}
				\left\langle y_{t+1}-y_{t}-\beta\widehat{\nabla}_{y} \tilde{f}_t\left(x_{t+1}, y_{t}\right),y-y_{t+1}\right\rangle\geq 0,\quad \forall y\in\mathcal{Y},\label{opt_yt+1_b}\\
				\left\langle y_{t}-y_{t-1}-\beta\widehat{\nabla}_{y} \tilde{f}_{t-1}\left(x_{t}, y_{t-1}\right),y-y_{t}\right\rangle\geq 0,\quad \forall y\in\mathcal{Y}.\label{opt_yt_b}
			\end{align}
			By letting $y=y_{t}$ in \eqref{opt_yt+1_b} and $y=y_{t+1}$ in \eqref{opt_yt_b} respectively, and multiplying $-\frac{1}{\beta}$ on both sides of two inequalities and rearranging terms, we obtain
			\begin{align}
				\left\langle\frac{1}{\beta}(y_{t+1}-y_{t})-\widehat{\nabla}_{y} \tilde{f}_t\left(x_{t+1}, y_{t}\right),y_{t+1}-y_{t}\right\rangle\leq 0,\label{1thop_yt+1_b}\\
				\left\langle\widehat{\nabla}_{y} \tilde{f}_{t-1}\left(x_{t}, y_{t-1}\right)-\frac{1}{\beta}(y_{t}-y_{t-1}), y_{t+1}-y_{t}\right\rangle\leq 0. \label{1thop_yt_b}
			\end{align}
			Denote $v_{t+1}=y_{t+1}-y_{t}-(y_{t}-y_{t-1})$. By adding \eqref{1thop_yt+1_b} and \eqref{1thop_yt_b}, we get
			\begin{align}\label{El7_b}
				\frac{1}{\beta}\left\langle v_{t+1},y_{t+1}-y_{t}\right\rangle
				&\leq\left\langle\widehat{\nabla}_{y}\tilde{f}_t(x_{t+1},y_{t})
				-\widehat{\nabla}_{y}\tilde{f}_{t-1}(x_{t},y_{t-1}),y_{t+1}
				-y_{t}\right\rangle\notag\\
				&=\left\langle\widehat{\nabla}_{y}\tilde{f}_t(x_{t+1},y_{t})
				-\widehat{\nabla}_{y}\tilde{f}_{t-1}(x_{t},y_{t}),y_{t+1}
				-y_{t}\right\rangle\notag\\
				&\quad +\left\langle\widehat{\nabla}_{y}\tilde{f}_{t-1}(x_{t},y_{t})
				-\widehat{\nabla}_{y}\tilde{f}_{t-1}(x_{t},y_{t-1}),
				y_{t+1}-y_{t}\right\rangle.
			\end{align}
			On the one hand,
			\begin{align}\label{El7_1_b}
				&\left\langle\widehat{\nabla}_{y}\tilde{f}_t(x_{t+1},y_{t})
				-\widehat{\nabla}_{y}\tilde{f}_{t-1}(x_{t},y_{t}),
				y_{t+1}-y_{t}\right\rangle\notag\\
				=&\left\langle\widehat{\nabla}_{y}f(x_{t+1},y_{t})
				-\lambda_{t}y_{t}
				-(\widehat{\nabla}_{y}f(x_{t},y_{t})
				-\lambda_{t-1}y_{t}),
				y_{t+1}-y_{t}\right\rangle\notag\\
				=&\left\langle\widehat{\nabla}_{y}f(x_{t+1},y_{t})
				-\widehat{\nabla}_{y}f(x_{t},y_{t}),
				y_{t+1}-y_{t}\right\rangle
				+(\lambda_{t-1}-\lambda_{t})\left\langle y_{t},
				y_{t+1}-y_{t}\right\rangle,
			\end{align}
			We estimate the first term in the r.h.s of \eqref{El7_1_b} as follows. By using the Cauchy-Schwarz inequality, Lemma \ref{varboundlemma_b}, and Assumption \ref{a4}, we have
			\begin{align}\label{El7_1_1_b}
				&\left\langle\widehat{\nabla}_{y}f(x_{t+1},y_{t})
				-\widehat{\nabla}_{y}f(x_{t},y_{t}),y_{t+1}-y_{t}\right\rangle\notag\\
				\leq&\frac{12}{\lambda_{t}}\|\widehat{\nabla}_{y}f(x_{t+1},y_{t})
				-\nabla_{y}f(x_{t+1},y_{t})\|^{2}
				+\frac{\lambda_{t}}{48}\|y_{t+1}-y_{t}\|^{2}\notag\\
				& +\frac{12}{\lambda_{t}}\|\nabla_{y}f(x_{t+1},y_{t})
				-\nabla_{y}f(x_{t},y_{t})\|^{2}
				+\frac{\lambda_{t}}{48}\|y_{t+1}-y_{t}\|^{2}\notag\\
				& +\frac{12}{\lambda_{t}}\|\nabla_{y}f(x_{t},y_{t})
				-\widehat{\nabla}_{y}f(x_{t},y_{t})\|^{2}
				+\frac{\lambda_{t}}{48}\|y_{t+1}-y_{t}\|^{2}\notag\\
				\leq&\frac{\lambda_{t}}{16}\|y_{t+1}-y_{t}\|^{2}
				+\frac{6d_yL_y^2\mu_2^2}{\lambda_{t}}
				+\frac{12}{\lambda_{t}}L_{y}^{2}\|x_{t+1}-x_{t}\|^{2}.
			\end{align}
			Moreover, it can be easily checked that
			\begin{equation}\label{El7_1_2_b}
				\left\langle y_{t},y_{t+1}-y_{t}\right\rangle=\frac{1}{2}\|y_{t+1}\|^{2}
				-\frac{1}{2}\|y_{t}\|^{2}-\frac{1}{2}\|y_{t+1}-y_{t}\|^{2}.
			\end{equation}
			Plugging \eqref{El7_1_1_b} and \eqref{El7_1_2_b} into \eqref{El7_1_b}, by Assumption \ref{ass:2}, we get
			\begin{align}\label{El7_1_new_b}
				&\left\langle\widehat{\nabla}_{y}\tilde{f}_t(x_{t+1},y_{t})
				-\widehat{\nabla}_{y}\tilde{f}_{t-1}(x_{t},y_{t}),
				y_{t+1}-y_{t}\right\rangle\notag\\
				\leq&\frac{\lambda_{t}}{16}\|y_{t+1}-y_{t}\|^{2}
				+\frac{6d_yL_y^2\mu_2^2}{\lambda_{t}}
				+\frac{12}{\lambda_{t}}L_{y}^{2}\|x_{t+1}-x_{t}\|^{2}+\frac{\lambda_{t-1}-\lambda_{t}}{2}\|y_{t+1}\|^{2}\nonumber\\
				&-\frac{\lambda_{t-1}-\lambda_{t}}{2}\|y_{t}\|^{2}.
			\end{align}
			Next, we estimate the second term in the r.h.s. of \eqref{El7_b} as follows. By the definition of $v_{t+1}$, we obtain
			\begin{align}\label{El7_2_b}
				&\left\langle\widehat{\nabla}_{y}\tilde{f}_{t-1}(x_{t},y_{t})
				-\widehat{\nabla}_{y}\tilde{f}_{t-1}(x_{t},y_{t-1}),y_{t+1}-y_{t}\right\rangle\notag\\
				=&\left\langle\widehat{\nabla}_{y}\tilde{f}_{t-1}(x_{t},y_{t})
				-\widehat{\nabla}_{y}\tilde{f}_{t-1}(x_{t},y_{t-1}),v_{t+1}\right\rangle\nonumber\\
				&+\left\langle\widehat{\nabla}_{y}\tilde{f}_{t-1}(x_{t},y_{t})
				-\widehat{\nabla}_{y}\tilde{f}_{t-1}(x_{t},y_{t-1}),y_{t}-y_{t-1}\right\rangle.
			\end{align}
			We first estimate the first term in the right side of \eqref{El7_2_b} as follows. By the Cauchy-Schwarz inequality and Lemma \ref{varboundlemma_b}, we have
			\begin{align}\label{El7_2_1_b}
				&\left\langle\widehat{\nabla}_{y}\tilde{f}_{t-1}(x_{t},y_{t})
				-\widehat{\nabla}_{y}\tilde{f}_{t-1}(x_{t},y_{t-1}),v_{t+1}\right\rangle\notag\\
				\leq& 2\beta\|\widehat{\nabla}_{y}\tilde{f}_{t-1}(x_{t},
				y_{t})-\nabla_{y}\tilde{f}_{t-1}(x_{t},y_{t})\|^{2}
				+\frac{1}{8\beta}\|v_{t+1}\|^{2}\nonumber\\
				&+\beta \|\nabla_{y}\tilde{f}_{t-1}(x_{t},y_{t})
				-\nabla_{y}\tilde{f}_{t-1}(x_{t},y_{t-1})\|^{2}
				+\frac{1}{4\beta}\|v_{t+1}\|^{2}\notag\\
				&
				+2\beta\|\nabla_{y}\tilde{f}_{t-1}(x_{t},
				y_{t-1})-\widehat{\nabla}_{y}\tilde{f}_{t-1}(x_{t},y_{t-1})\|^{2}
				+\frac{1}{8\beta}\|v_{t+1}\|^{2}\notag\\
				\leq& \beta d_{y}L_{y}^{2}\mu_{2}^{2}
				+\beta\|\nabla_{y}\tilde{f}_{t-1}(x_{t},y_{t})
				-\nabla_{y}\tilde{f}_{t-1}(x_{t},y_{t-1})\|^{2}+\frac{1}{2\beta}\|v_{t+1}\|^{2}.
			\end{align}
			On the other hand, we can rewrite the second term in the right side of \eqref{El7_2_b} as
			\begin{align}\label{El7_fiveterms_b}
				&\left\langle\widehat{\nabla}_{y}\tilde{f}_{t-1}(x_{t},y_{t})
				-\widehat{\nabla}_{y}\tilde{f}_{t-1}(x_{t},y_{t-1}),y_{t}-y_{t-1}\right\rangle\notag\\
				=&\left\langle\widehat{\nabla}_{y}\tilde{f}_{t-1}(x_{t},y_{t})
				-\nabla_{y}\tilde{f}_{t-1}(x_{t},y_{t}),y_{t}-y_{t-1}\right\rangle\nonumber\\
				&+\left\langle\nabla_{y}\tilde{f}_{t-1}(x_{t},y_{t})
				-\nabla_{y}\tilde{f}_{t-1}(x_{t},y_{t-1}),y_{t}-y_{t-1}\right\rangle\notag\\
				&+\left\langle\nabla_{y}\tilde{f}_{t-1}(x_{t},y_{t-1})
				-\widehat{\nabla}_{y}\tilde{f}_{t-1}(x_{t},y_{t-1}),y_{t}-y_{t-1}\right\rangle.
			\end{align}
			By the Cauchy-Schwarz inequality and Lemma \ref{varboundlemma_b}, we obtain
			\begin{align}\label{El7_2_2_2_b}
				&\left\langle\widehat{\nabla}_{y}\tilde{f}_{t-1}(x_{t},y_{t})
				-\nabla_{y}\tilde{f}_{t-1}(x_{t},y_{t}),y_{t}-y_{t-1}\right\rangle\nonumber\\
				\leq&\frac{4}{\lambda_{t}}\|\widehat{\nabla}_{y}\tilde{f}_{t-1}(x_{t},y_{t})
				-\nabla_{y}\tilde{f}_{t-1}(x_{t},y_{t})\|^{2}+\frac{\lambda_{t}}{16}
				\|y_{t}-y_{t-1}\|^{2}\notag\\
				\leq &\frac{d_{y}L_{y}^{2}\mu_{2}^{2}}{\lambda_{t}}
				+\frac{\lambda_{t}}{16}\|y_{t}-y_{t-1}\|^{2}.
			\end{align}
			Noticing that $\nabla_{y}\tilde{f}_{t-1}(x_{t},y)=
			\nabla_{y}f(x_{t},y)-\lambda_{t-1}y$ and by Assumption \ref{a4}, we have 
			that $\tilde{f}_{t-1}(x,y)$ has Lipschitz continuous gradients with $\tilde{L}_{y}$, where $\tilde{L}_{y}=L_{y}+\lambda_{1}$. 
			In view of the strong concavity of $\tilde{f}(x,y)$ w.r.t. $y$, from \cite[Theorem 2.1.12]{nesterov2018lectures} we obtain
			\begin{align}
				&\left\langle\nabla_{y}\tilde{f}_{t-1}(x_{t},y_{t})
				-\nabla_{y}\tilde{f}_{t-1}(x_{t},y_{t-1}),
				y_{t}-y_{t-1}\right\rangle\nonumber\\
				\leq&-\frac{1}{\tilde{L}_{y}+\lambda_{t-1}}\|\nabla_{y}\tilde{f}_{t-1}(x_{t},y_{t})
				-\nabla_{y}\tilde{f}_{t-1}(x_{t},y_{t-1})\|^{2}-\frac{\tilde{L}_{y}\lambda_{t-1}}{\tilde{L}_{y}+\lambda_{t-1}}
				\|y_{t}-y_{t-1}\|^{2}.\label{El7_2_2_3_b}
			\end{align}
			Similar to \eqref{El7_2_2_2_b}, the third term in the right side of \eqref{El7_fiveterms_b} can be estimated as
			\begin{align}
				\left\langle\nabla_{y}\tilde{f}_{t-1}(x_{t},y_{t-1})
				-\widehat{\nabla}_{y}\tilde{f}_{t-1}(x_{t},y_{t-1}),
				y_{t}-y_{t-1}\right\rangle
				\leq\frac{d_{y}L_{y}^{2}\mu_{2}^{2}}{\lambda_{t}}
				+\frac{\lambda_{t}}{16}\|y_{t}-y_{t-1}\|^{2}.\label{El7_2_2_4_b}
			\end{align}
			Next, by plugging \eqref{El7_2_2_2_b} - \eqref{El7_2_2_4_b} into \eqref{El7_fiveterms_b}, we obtain
			\begin{align}\label{El7_fiveterms_new_b}
				&\left\langle\widehat{\nabla}_{y}\tilde{f}_{t-1}(x_{t},y_{t})
				-\widehat{\nabla}_{y}\tilde{f}_{t-1}(x_{t},y_{t-1}),y_{t}-y_{t-1}\right\rangle\notag\\
				\leq& 
				-\left(\frac{\tilde{L}_{y}\lambda_{t-1}}{\tilde{L}_{y}+\lambda_{t-1}}-\frac{\lambda_{t}}{8}\right)\|y_{t}-y_{t-1}\|^{2}\nonumber\\
				&-\frac{1}{\tilde{L}_{y}+\lambda_{t-1}}\|\nabla_{y}\tilde{f}_{t-1}(x_{t},y_{t})
				-\nabla_{y}\tilde{f}_{t-1}(x_{t},y_{t-1})\|^{2}+\frac{2d_{y}L_{y}^{2}\mu_{2}^{2}}{\lambda_{t}}.
			\end{align}
			Plugging \eqref{El7_2_1_b} and \eqref{El7_fiveterms_new_b} into \eqref{El7_2_b} and rearranging terms, we get that
			\begin{align}\label{El7_2_new_b}
				&\left\langle\widehat{\nabla}_{y}\tilde{f}_{t-1}(x_{t},y_{t})
				-\widehat{\nabla}_{y}\tilde{f}_{t-1}(x_{t},y_{t-1}),y_{t+1}-y_{t}\right\rangle\notag\\
				\leq&\frac{1}{2\beta}\|v_{t+1}\|^{2} +\left(\beta-\frac{1}{\tilde{L}_{y}+\lambda_{t-1}}\right) \|\nabla_{y}\tilde{f}_{t-1}(x_{t},
				y_{t})-\nabla_{y}\tilde{f}_{t-1}(x_{t},y_{t-1})\|^{2}\notag\\
				& -\left(\frac{\tilde{L}_{y}\lambda_{t-1}}{\tilde{L}_{y}+\lambda_{t-1}}-\frac{\lambda_{t}}{8}\right) \|y_{t}-y_{t-1}\|^{2} +\left(\beta+\frac{2}{\lambda_{t}}\right)d_{y}L_{y}^{2}\mu_{2}^{2}\notag\\
				\leq&\frac{1}{2\beta}\|v_{t+1}\|^{2}
				-\frac{3\lambda_{t}}{8}\|y_{t}-y_{t-1}\|^{2}
				+\left(\beta+\frac{2}{\lambda_{t}}\right)d_{y}L_{y}^{2}\mu_{2}^{2},
			\end{align}
			where the last inequality is due to the choice of $\beta \leq \frac{1}{L_{y}+2\lambda_{1}} \leq \frac{1}{\tilde{L}_{y}+\lambda_{t-1}}$ and 
			$\frac{\tilde{L}_{y}\lambda_{t-1}}{\tilde{L}_{y}+\lambda_{t-1}}\geq
			\frac{\tilde{L}_{y}\lambda_{t-1}}{2\tilde{L}_{y}}=\frac{\lambda_{t-1}}{2}\ge\frac{\lambda_{t}}{2}.$
			Then plugging \eqref{El7_1_new_b} and \eqref{El7_2_new_b} into \eqref{El7_b} and rearranging terms, we obtain
			\begin{align}\label{El7_n_b}
				&\frac{1}{\beta}\left\langle v_{t+1},y_{t+1}-y_{t}\right\rangle\notag\\
				\leq& \frac{1}{2\beta}\|v_{t+1}\|^{2} -\frac{3\lambda_{t}}{8}\|y_{t}-y_{t-1}\|^{2}
				+\frac{\lambda_{t}}{16}\|y_{t+1}-y_{t}\|^{2}
				+\frac{\lambda_{t-1}-\lambda_{t}}{2}\|y_{t+1}\|^{2}
				\notag\\
				& -\frac{\lambda_{t-1}-\lambda_{t}}{2}\|y_{t}\|^{2}+\frac{12}{\lambda_{t}}L_{y}^{2}\|x_{t+1}-x_{t}\|^{2} +\left(\beta+\frac{8}{\lambda_{t}}\right)d_{y}L_{y}^{2}\mu_{2}^{2}.
			\end{align}
			It is easy to check that
			\begin{equation}\label{threepoint_b}
				\frac{1}{\beta}\left\langle v_{t+1},y_{t+1}-y_{t}\right\rangle
				=\frac{1}{2\beta}\|v_{t+1}\|^{2}
				+\frac{1}{2\beta}\|y_{t+1}-y_{t}\|^{2}
				-\frac{1}{2\beta}\|y_{t}-y_{t-1}\|^{2}.
			\end{equation}
			By plugging \eqref{threepoint_b} into \eqref{El7_n_b} and rearranging terms, we can get
			\begin{align*}
				&\frac{1}{2\beta}\|y_{t+1}-y_{t}\|^{2}
				-\frac{\lambda_{t-1}-\lambda_{t}}{2}\|y_{t+1}\|^{2}\\
				\leq&\frac{1}{2\beta}\|y_{t}-y_{t-1}\|^{2}
				-\frac{\lambda_{t-1}-\lambda_{t}}{2}\|y_{t}\|^{2}
				-\frac{3\lambda_{t}}{8}\|y_{t}-y_{t-1}\|^{2}
				+\frac{\lambda_{t}}{16}\|y_{t+1}-y_{t}\|^{2}
				\\
				&
				+\frac{12}{\lambda_{t}}L_{y}^{2}\|x_{t+1}-x_{t}\|^{2}+\left(\beta+\frac{8}{\lambda_{t}}\right)d_{y}L_{y}^{2}\mu_{2}^{2}.
			\end{align*}
			By multiplying both sides of the inequality by $\frac{32}{\beta\lambda_{t}}$, we can get
			\begin{align}\label{l7_ori_b}
				&\frac{16}{\beta^{2}\lambda_{t}}\left\|y_{t+1}-y_{t}\right\|^{2}
				-\frac{16}{\beta}\left(\frac{\lambda_{t-1}}{\lambda_{t}}-1\right)\left\|y_{t+1}\right\|^{2}\notag\\
				\leq&\frac{16}{\beta^{2}\lambda_{t}}\left\|y_{t}-y_{t-1}\right\|^{2}
				-\frac{16}{\beta}\left(\frac{\lambda_{t-1}}{\lambda_{t}}-1\right)\left\|y_{t}\right\|^{2}-\frac{12}{\beta}\left\|y_{t}-y_{t-1}\right\|^{2}\notag\\
				&+\frac{2}{\beta}\left\|y_{t+1}-y_{t}\right\|^{2} +\frac{384}{\beta\lambda_{t}^{2}}L_{y}^{2}\left\|x_{t+1}-x_{t}\right\|^{2}
				+\left(\frac{32}{\lambda_{t}}+\frac{256}{\beta\lambda_{t}^{2}}\right)d_{y}L_{y}^{2}\mu_{2}^{2}.
			\end{align}
			Finally, by adding two terms, i.e., $\frac{16}{\beta^{2}\lambda_{t+1}}\left\|y_{t+1}-y_{t}\right\|^{2}$ and $\frac{16}{\beta}\left(\frac{\lambda_{t-2}}{\lambda_{t-1}}-1\right)\left\|y_{t}\right\|^{2}$ to both sides of \eqref{l7_ori_b}, subsequently rearranging terms, the proof is completed.	
		\end{proof}
		
		\begin{lemma}\label{decentxy_block}
			Suppose that Assumption \ref{a4} holds. Assume $\beta \leq \frac{1}{L_{y}+2\lambda_{1}}$. Then, for the sequence $\{(x_{t},y_{t})\}$ generated by Algorithm \ref{algo1}, and for any $t\geq 3$, we have
			\begin{align*}
					\mathcal{L}(x_{t+1},y_{t+1})\leq 	& 	\mathcal{L}(x_{t},y_{t})
				-c_{t}\|y_{t+1}-y_{t}\|^{2}
				-m_{t}\|x_{t+1}-x_{t}\|^{2}+\frac{\beta d_xL_x^2\mu_{1,t}^2}{8}\\
				&+ \frac{\lambda_{t-1}-\lambda_{t}}{2}\|y_{t+1}\|^{2}
				+\frac{16}{\rho}\left(\frac{\lambda_{t-2}}{\lambda_{t-1}}-\frac{\lambda_{t-1}}{\lambda_{t}}\right)
				\|y_{t}\|^{2}
				+ r_{t}d_{y}L_{y}^{2}\mu_{2}^{2},
			\end{align*}
			where $	\mathcal{L}(x_{t},y_{t})=f(x_{t},y_{t})+\frac{16}{\beta^{2}\lambda_{t}}\|y_{t}-y_{t-1}\|^{2}
			-\frac{16}{\beta}\left(\frac{\lambda_{t-2}}{\lambda_{t-1}}-1\right)\|y_{t}\|^{2} -\frac{\lambda_{t-1}}{2}\|y_{t}\|^{2}
			-\frac{11}{\beta}\|y_{t}-y_{t-1}\|^{2}$ and  the parameters are chosen as follows,
			\begin{align}
				c_{t} &= \frac{15}{2\beta}-\frac{16}{\beta^{2}}\left(\frac{1}{\lambda_{t+1}}
				-\frac{1}{\lambda_{t}}\right),\quad m_{t} = \frac{1}{\alpha_t}-\left(\frac{L_{x}}{2}+\frac{1}{2\beta}+\frac{384}{\beta\lambda_{t}^{2}}L_{y}^{2}+\frac{3\beta L_{y}^{2}}{2}\right),\label{c_m_b}\\
				r_{t} &= \frac{19\beta}{8}+\frac{32}{\lambda_{t}}+\frac{256}{\beta\lambda_{t}^{2}}.\label{p_r_b}
			\end{align}
		\end{lemma}
		\begin{proof}
			In view of the concavity of $\tilde{f}_t(x,y)$ on $y$, we have 
			\begin{equation}
				\tilde{f}_t(x_{t+1},y_{t+1})-\tilde{f}_t(x_{t+1},y_{t})\leq \left\langle\nabla_{y}\tilde{f}_t(x_{t+1},y_{t}),
				y_{t+1}-y_{t}\right\rangle.\label{add:2}
			\end{equation}
			By splitting the r.h.s. of \eqref{add:2} into three terms, and \eqref{1thop_yt_b}, we obtain that
			\begin{align}\label{l8_0_b}
				\tilde{f}_t(x_{t+1},y_{t+1})-\tilde{f}_t(x_{t+1},y_{t})
				\leq&\left\langle\nabla_{y}\tilde{f}_t(x_{t+1},y_{t})
				-\widehat{\nabla}_{y}\tilde{f}_t(x_{t+1},y_{t}),
				y_{t+1}-y_{t}\right\rangle\notag\\
				&
				+\left\langle\widehat{\nabla}_{y}\tilde{f}_t(x_{t+1},y_{t})
				-\widehat{\nabla}_{y}\tilde{f}_{t-1}(x_{t},y_{t-1}),
				y_{t+1}-y_{t}\right\rangle\notag\\
				&
				+\frac{1}{\beta}\left\langle y_{t}-y_{t-1},
				y_{t+1}-y_{t}\right\rangle.
			\end{align}
			By using $\tilde{f}_t(x_{t+1},y_{t+1})
			=f(x_{t+1},y_{t+1})-\frac{\lambda_{t}}{2}\|y_{t+1}\|^{2}$ and $\tilde{f}_t(x_{t+1},y_{t})
			=f(x_{t+1},y_{t})-\frac{\lambda_{t}}{2}\|y_{t}\|^{2}$, we obtain
			\begin{align}\label{El8_0_new_b}
				&f(x_{t+1},y_{t+1})
				-\frac{\lambda_{t}}{2}\|y_{t+1}\|^{2}
				-\left(f(x_{t+1},y_{t})
				-\frac{\lambda_{t}}{2}\|y_{t}\|^{2}\right)\notag\\
				\leq&\left\langle\nabla_{y}\tilde{f}_t(x_{t+1},y_{t})
				-\widehat{\nabla}_{y}\tilde{f}_t(x_{t+1},y_{t}),
				y_{t+1}-y_{t}\right\rangle\notag\\
				& +\left\langle\widehat{\nabla}_{y}\tilde{f}_t(x_{t+1},y_{t})
				-\widehat{\nabla}_{y}\tilde{f}_{t-1}(x_{t},y_{t-1}),
				y_{t+1}-y_{t}\right\rangle
				+\frac{1}{\beta}\left\langle y_{t}-y_{t-1},
				y_{t+1}-y_{t}\right\rangle.
			\end{align}
			Next, we estimate the upper bound of the r.h.s. of \eqref{El8_0_new_b}. Firstly, by the Cauchy-Schwarz inequality and Lemma \ref{varboundlemma_b}, we have
			\begin{align}\label{El8_1_1_b}
				&\left\langle\nabla_{y}\tilde{f}_t(x_{t+1},y_{t})
				-\widehat{\nabla}_{y}\tilde{f}_t(x_{t+1},y_{t}),
				y_{t+1}-y_{t}\right\rangle\notag\\
				\leq&\frac{\beta}{2}\|\nabla_{y}f(x_{t+1},y_{t})
				-\widehat{\nabla}_{y}f(x_{t+1},y_{t})\|^{2}
				+\frac{1}{2\beta}\|y_{t+1}-y_{t}\|^{2}\notag\\
				\leq&\frac{\beta d_{y}L_{y}^{2}\mu_{2}^{2}}{8} +\frac{1}{2\beta}\|y_{t+1}-y_{t}\|^{2}.
			\end{align}
			Secondly, 
			\begin{align}\label{El8_1_2_b}
				&\left\langle\widehat{\nabla}_{y}\tilde{f}_t(x_{t+1},y_{t})
				-\widehat{\nabla}_{y}\tilde{f}_{t-1}(x_{t},y_{t-1}),
				y_{t+1}-y_{t}\right\rangle\notag\\
				=&\left\langle\widehat{\nabla}_{y}\tilde{f}_t(x_{t+1},y_{t})
				-\widehat{\nabla}_{y}\tilde{f}_{t-1}(x_{t},y_{t}),
				y_{t+1}-y_{t}\right\rangle\nonumber\\ &+\left\langle\widehat{\nabla}_{y}\tilde{f}_{t-1}(x_{t},y_{t})
				-\widehat{\nabla}_{y}\tilde{f}_{t-1}(x_{t},y_{t-1}),
				v_{t+1}\right\rangle\notag\\
				& +\left\langle\widehat{\nabla}_{y}\tilde{f}_{t-1}(x_{t},y_{t})
				-\widehat{\nabla}_{y}\tilde{f}_{t-1}(x_{t},y_{t-1}),
				y_{t}-y_{t-1}\right\rangle.
			\end{align}
			Next, by the Cauchy-Schwarz inequality and Assumption \ref{a4}, using Lemma \ref{varboundlemma_b}, we obtain
			\begin{align}\label{El8_1_2_1_1_b}
				&\left\langle\widehat{\nabla}_{y}f(x_{t+1},y_{t})
				-\widehat{\nabla}_{y}f(x_{t},y_{t}),
				y_{t+1}-y_{t}\right\rangle\notag\\
				\leq&\frac{3\beta}{2}\|\widehat{\nabla}_{y}f(x_{t+1},y_{t})
				-\nabla_{y}f(x_{t+1},y_{t})\|^{2}
				+\frac{1}{6\beta}\|y_{t+1}-y_{t}\|^{2}\notag\\
				& +\frac{3\beta}{2}\|\nabla_{y}f(x_{t+1},y_{t})
				-\nabla_{y}f(x_{t},y_{t})\|^{2}
				+\frac{1}{6\beta}\|y_{t+1}-y_{t}\|^{2}\notag\\
				& +\frac{3\beta}{2}\|\nabla_{y}f(x_{t},y_{t})
				-\widehat{\nabla}_{y}f(x_{t},y_{t})\|^{2}
				+\frac{1}{6\beta}\|y_{t+1}-y_{t}\|^{2}\notag\\
				\leq& \frac{3\beta d_yL_y^2\mu_2^2}{4}
				+\frac{3\beta L_{y}^{2}}{2}\|x_{t+1}-x_{t}\|^{2} +\frac{1}{2\beta}\|y_{t+1}-y_{t}\|^{2}.
			\end{align}
			By plugging \eqref{El8_1_2_1_1_b} and \eqref{El7_1_2_b} into \eqref{El7_1_b} and Assumption \ref{ass:2}, we obtain
			\begin{align}\label{El8_1_2_1_new_b}
				&\left\langle\widehat{\nabla}_{y}\tilde{f}_t(x_{t+1},y_{t})
				-\widehat{\nabla}_{y}\tilde{f}_{t-1}(x_{t},y_{t}),
				y_{t+1}-y_{t}\right\rangle\notag\\
				\leq& \frac{3\beta d_yL_y^2\mu_2^2}{4}
				+\frac{3\beta L_{y}^{2}}{2}\|x_{t+1}-x_{t}\|^{2} +\frac{1}{2\beta}\|y_{t+1}-y_{t}\|^{2}\nonumber\\
				&+\frac{\lambda_{t-1}-\lambda_{t}}{2}(\|y_{t+1}\|^{2}
				-\|y_{t}\|^{2}).
			\end{align}
			Next, the third term in the r.h.s. of \eqref{El8_1_2_b} can be rewritten the same as in \eqref{El7_fiveterms_b}, which has three terms in the right side. For the first and third term, similar to the proof of \eqref{El7_2_2_2_b}, we can get
			\begin{align}\label{El8_2_2_2_b}
				&\left\langle\widehat{\nabla}_{y}\tilde{f}_{t-1}(x_{t},y_{t})
				-\nabla_{y}\tilde{f}_{t-1}(x_{t},y_{t}),y_{t}-y_{t-1}\right\rangle\leq\frac{\beta d_yL_y^2\mu_2^2}{4}
				+\frac{1}{4\beta}\|y_{t}-y_{t-1}\|^{2}
			\end{align}
			and
			\begin{align}\label{El8_2_2_4_b}
				\left\langle\nabla_{y}\tilde{f}_{t-1}(x_{t},y_{t-1})
				-\widehat{\nabla}_{y}\tilde{f}_{t-1}(x_{t},y_{t-1}),
				y_{t}-y_{t-1}\right\rangle
				\leq\frac{\beta d_yL_y^2\mu_2^2}{4}
				+\frac{1}{4\beta}\|y_{t}-y_{t-1}\|^{2}.
			\end{align}
			By plugging \eqref{El8_2_2_2_b}, \eqref{El7_2_2_3_b} and \eqref{El8_2_2_4_b} into \eqref{El7_fiveterms_b}, we obtain
			\begin{align}\label{five_terms_new_b}
				&\left\langle\widehat{\nabla}_{y}\tilde{f}_{t-1}(x_{t},y_{t})
				-\widehat{\nabla}_{y}\tilde{f}_{t-1}(x_{t},y_{t-1}),y_{t}-y_{t-1}\right\rangle\notag\\
				\leq& 
				\frac{1}{2\beta}\|y_{t}-y_{t-1}\|^{2}
				-\frac{1}{\tilde{L}_{y}+\lambda_{t-1}} \|\nabla_{y}\tilde{f}_{t-1}(x_{t},y_{t})
				-\nabla_{y}\tilde{f}_{t-1}(x_{t},y_{t-1})\|^{2}\nonumber\\
				&+\frac{\beta d_yL_y^2\mu_2^2}{2}.
			\end{align}
			By plugging \eqref{El8_1_2_1_new_b}, \eqref{El7_2_1_b} and \eqref{five_terms_new_b} into \eqref{El8_1_2_b} and rearranging terms, using $\beta \leq \frac{1}{L_{y}+2\lambda_{1}}$, we can get
			\begin{align}\label{El8_1_2_new_b}
				&\left\langle\widehat{\nabla}_{y}\tilde{f}_t(x_{t+1},y_{t}) 
				-\widehat{\nabla}_{y}\tilde{f}_{t-1}(x_{t},y_{t-1}),
				y_{t+1}-y_{t}\right\rangle\notag\\
				\leq&\frac{1}{2\beta}\|y_{t+1}-y_{t}\|^{2}
				+\frac{1}{2\beta}\|y_{t}-y_{t-1}\|^{2}
				+\frac{3\beta L_{y}^{2}}{2}\|x_{t+1}-x_{t}\|^{2}+\frac{1}{2\beta}\|v_{t+1}\|^{2}\notag\\
				&
				+\frac{\lambda_{t-1}-\lambda_{t}}{2}(\|y_{t+1}\|^{2}
				-\|y_{t}\|^{2})
				+ \frac{9\beta d_yL_y^2\mu_2^2}{4}.
			\end{align}
			Next, by plugging \eqref{El8_1_1_b}, \eqref{El8_1_2_new_b} into \eqref{El8_0_new_b}, we obtain
			\begin{align}\label{EL81_0_b}
				&f(x_{t+1},y_{t+1})
				-\frac{\lambda_{t}}{2}\|y_{t+1}\|^{2}
				-\left(f(x_{t+1},y_{t})
				-\frac{\lambda_{t}}{2}\|y_{t}\|^{2}\right)\notag\\
				\leq&\frac{1}{\beta}\|y_{t+1}-y_{t}\|^{2}
				+\frac{1}{2\beta}\|y_{t}-y_{t-1}\|^{2}
				+\frac{3\beta L_{y}^{2}}{2}\|x_{t+1}-x_{t}\|^{2}+\frac{1}{2\beta}\|v_{t+1}\|^{2}\notag\\
				&
				+\frac{\lambda_{t-1}-\lambda_{t}}{2}(\|y_{t+1}\|^{2}
				-\|y_{t}\|^{2}) +\frac{19\beta d_yL_y^2\mu_2^2}{8}
				+\frac{1}{\beta}\left\langle y_{t}-y_{t-1}, y_{t+1}-y_{t}\right\rangle.
			\end{align}
			By plugging the following equality into the above inequality, i.e.,
			\begin{equation}
				\left\langle y_{t}-y_{t-1},
				y_{t+1}-y_{t}\right\rangle=
				\frac{1}{2}\|y_{t+1}-y_{t}\|^{2}
				+\frac{1}{2}\|y_{t}-y_{t-1}\|^{2}
				-\frac{1}{2}\|v_{t+1}\|^{2},
			\end{equation}
			we get that
			\begin{align}\label{EL81_b}
				&f(x_{t+1},y_{t+1})
				-\frac{\lambda_{t}}{2}\|y_{t+1}\|^{2}
				-\left(f(x_{t+1},y_{t})
				-\frac{\lambda_{t}}{2}\|y_{t}\|^{2}\right)\notag\\
				\leq &\frac{3}{2\beta}\|y_{t+1}-y_{t}\|^{2}
				+\frac{1}{\beta}\|y_{t}-y_{t-1}\|^{2}
				+\frac{3\beta L_{y}^{2}}{2}\|x_{t+1}-x_{t}\|^{2}\notag\\
				&+\frac{\lambda_{t-1}-\lambda_{t}}{2}(\|y_{t+1}\|^{2}
				-\|y_{t}\|^{2}) 
				+\frac{19\beta d_yL_y^2\mu_2^2}{8}.
			\end{align}
			By adding \eqref{disy_block} and \eqref{EL81_b}, and rearranging terms, we have
			\begin{align}\label{l8_fina_b}
				&f(x_{t+1},y_{t+1})
				-\frac{\lambda_{t}}{2}\|y_{t+1}\|^{2}
				+\frac{16}{\beta^{2}\lambda_{t+1}}\|y_{t+1}-y_{t}\|^{2}-\frac{11}{\beta}\|y_{t+1}-y_{t}\|^{2}\notag\\
				&  -\frac{16}{\beta}\left(\frac{\lambda_{t-1}}{\lambda_{t}}-1\right)\|y_{t+1}\|^{2}\notag\\
				\leq&f(x_{t+1},y_{t})
				-\frac{\lambda_{t-1}}{2}\|y_{t}\|^{2}
				+\frac{16}{\beta^{2}\lambda_{t}}\|y_{t}-y_{t-1}\|^{2} -\frac{11}{\beta}\|y_{t}-y_{t-1}\|^{2} \notag\\
				& -\frac{16}{\beta}\left(\frac{\lambda_{t-2}}{\lambda_{t-1}}-1\right)\|y_{t}\|^{2} -\left[\frac{15}{2\beta}-\frac{16}{\beta^{2}}\left(\frac{1}{\lambda_{t+1}}
				-\frac{1}{\lambda_{t}}\right)\right]\|y_{t+1}-y_{t}\|^{2}
				\notag\\
				&+\left(\frac{384 L_{y}^{2}}{\beta\lambda_{t}^{2}}+\frac{3\beta L_{y}^{2}}{2}\right)
				\|x_{t+1}-x_{t}\|^{2} +\frac{\lambda_{t-1}-\lambda_{t}}{2}\|y_{t+1}\|^{2}
				\nonumber\\
				&+\frac{16}{\beta}\left(\frac{\lambda_{t-2}}{\lambda_{t-1}}-\frac{\lambda_{t-1}}{\lambda_{t}}\right)
				\|y_{t}\|^{2}+\left(\frac{19\beta}{8}+\frac{32}{\lambda_{t}}+\frac{256}{\beta\lambda_{t}^{2}}\right)d_{y}L_{y}^{2}\mu_{2}^{2}.
			\end{align}
			Then adding \eqref{decentx_block} to both sides of the above inequality, we complete the proof.	
		\end{proof}

		Now, we are ready to give the main theorem in this section.
		\begin{theorem}\label{thm:main_b}
			Suppose that Assumption \ref{a4} holds. The sequence $\{(x_{t},y_{t})\}$ is generated by Algorithm \ref{algo1}. Let $\lambda_{t}=\frac{9}{20\beta t^{1/4}}, \quad \beta\leq\frac{1}{10L_{y}},\quad \frac{1}{\alpha_t}=\frac{L_{x}}{2}+\frac{1}{2\beta}+\frac{768}{\beta\lambda_{t}^{2}}L_{y}^{2}+\frac{3\beta L_{y}^{2}}{2} .$
			Let $\theta =16+\frac{9^{4}\left[4(\frac{L_{x}}{2}+\frac{1}{2\beta}+\frac{3\beta L_{y}^{2}}{2} )^{2} +3L_{y}^{2}\right]}{\left(384\times 20^{2}\beta L_{y}^{2}\right)^{2}}$,
			$s =\max\left\{\theta, \frac{27^{2}}{448\times 20^{2}\beta^{2}L_{y}^{2}}\right\}$, $C_{1}=\frac{27}{128\times20^{2}\beta L_{y}^{2}}$, $C_2=\max\{\frac{7\beta}{54}, C_1\}$.
			Then for any $t\geq 3$ and for any given $\varepsilon \in (0,1)$, if we set
			\begin{align}
				\mu_{1,t}=\frac{\varepsilon}{L_xt^{1/4}}\sqrt{\frac{C_1}{2s(4C_2+\beta)d_x}},\quad
				\mu_{2}=\frac{\varepsilon}{4L_{y}\sqrt{d_y}}\sqrt{\frac{1}{\frac{10837s}{81L_yC_1}+\frac{3}{4}}},\label{mu2_t_b}
			\end{align}
			then
			$$T(\varepsilon) \leq \max\left\{\left(\frac{8sC}{\varepsilon^{2}C_{1}}+2\right)^{2}, \frac{ 9^{4}\sigma_{y}^{4}}{10^{4}\beta^{4}\varepsilon^{4}}\right\},$$
			where $C=\bar{f}-\underline{f}+ \frac{149}{\beta}\sigma_{y}^{2} $ and $\sigma_{y}=\max\{\|y\||y\in\mathcal{Y}\}$.
		\end{theorem}
		
		\begin{proof}
			Form the setting of $\lambda_{t}$ and $\beta$, it holds that
			$\frac{1}{\lambda_{t+1}}-\frac{1}{\lambda_{t}}\leq\frac{4\beta}{9}$ and  $\frac{1}{\alpha_t}>\frac{L_{x}}{2}+\frac{1}{2\beta}+\frac{384}{\beta\lambda_{t}^{2}}L_{y}^{2}+\frac{3\beta L_{y}^{2}}{2},$
			and $m_{t}=\frac{384}{\beta\lambda_{t}^{2}}L_{y}^{2}$. Then we have $c_{t}\geq\frac{7}{18\beta}$ which is defined as in \eqref{c_m_b}.  Note that $r_{t}>0$ defined as in \eqref{p_r_b}. On the one hand, by Lemma \ref{decentxy_block}, we obtain
			\begin{align}\label{potient_block}
					\mathcal{L}(x_{t+1},y_{t+1})\leq& 	\mathcal{L}(x_{t},y_{t})
				-\frac{7}{18\beta}\|y_{t+1}-y_{t}\|^{2}
				-m_{t}\|x_{t+1}-x_{t}\|^{2}+\frac{\beta d_xL_x^2\mu_{1,t}^2}{8}\notag\\
				&
				+\frac{\lambda_{t-1}-\lambda_{t}}{2}\sigma_{y}^{2}
				+\frac{16}{\beta}\left(\frac{\lambda_{t-2}}{\lambda_{t-1}}-\frac{\lambda_{t-1}}{\lambda_{t}}\right)
				\sigma_{y}^{2}+r_{t}d_{y}L_{y}^{2}\mu_{2}^{2},
			\end{align}
			where $\frac{\lambda_{t-1}-\lambda_{t}}{2} >0$ and $\frac{16}{\beta}\left(\frac{\lambda_{t-2}}{\lambda_{t-1}}-\frac{\lambda_{t-1}}{\lambda_{t}}\right) >0$.
			On the other hand, by nonexpansiveness of the projection operator and using the Cauchy-Schwarz inequality, we have
			\begin{align}\label{xGk_block}
				\left\|(\nabla\tilde{\mathbf{G}}_{t})_{x}\right\|
				&=\frac{1}{\alpha_t}\left\|x_{t+1}-x_{t} - (x_{t+1}-\operatorname{Proj}_{\mathcal{X}}(x_{t}-\alpha_t\nabla_{x}
				f(x_{t},y_{t})))\right\|\notag\\
				&\leq\frac{1}{\alpha_t}\left\|x_{t+1}-x_{t}\right\|
				+\left\|\widehat{\nabla}_x \tilde{f}_t(x_{t},y_{t})
				-\nabla_{x}f(x_{t},y_{t})\right\|,
			\end{align}
			and
			\begin{align}\label{yG_b}
				\left\|(\nabla\tilde{\mathbf{G}}_{t})_{y}\right\|
				&\leq\frac{1}{\beta}\left\|y_{t+1}-y_{t}\right\|
				+\frac{1}{\beta}\left\|y_{t+1}-\operatorname{proj}_{\mathcal{Y}}(y_{t}+\beta\nabla_{y}
				\tilde{f}_t(x_{t},y_{t}))\right\|\notag\\
				&\leq\frac{1}{\beta}\left\|y_{t+1}-y_{t}\right\|
				+\left\|\widehat{\nabla}_{y}f(x_{t+1},y_{t})-\lambda_{t}y_{t}
				-(\nabla_{y}f(x_{t},y_{t})-\lambda_{t}y_{t})\right\|\notag\\
				&\leq\frac{1}{\beta}\left\|y_{t+1}-y_{t}\right\|
				+\left\|\widehat{\nabla}_{y}f(x_{t+1},y_{t})-\nabla_{y}
				f(x_{t},y_{t})\right\|\nonumber\\
				&\le \frac{1}{\beta}\left\|y_{t+1}-y_{t}\right\|
				+\left\|\widehat{\nabla}_{y}f(x_{t+1},y_{t})-\nabla_{y}
				f(x_{t+1},y_{t})\right\|\nonumber\\
				&+\left\|\nabla_{y}f(x_{t+1},y_{t})-\nabla_{y}
				f(x_{t},y_{t})\right\|.
			\end{align}    
			Next, by using the Cauchy-Schwarz inequality and Lemma \ref{varboundlemma_b}, we obtain
			\begin{align}\label{xGk_block_E}
				\left\|(\nabla\tilde{\mathbf{G}}_{t})_{x}\right\|^{2}
				\leq \frac{2}{\alpha_t^2}\left\|x_{t+1}-x_{t}\right\|^{2} + \frac{d_xL_x^2\mu_{1,t}^2}{2} .
			\end{align}
			Similarly, we have
			\begin{align}
				\left\|(\nabla\tilde{\mathbf{G}}_{t})_{y}\right\|^{2}
				\leq \frac{3}{\beta^{2}}\left\|y_{t+1}-y_{t}\right\|^{2}+3L_y^2\left\|x_{t+1}-x_{t}\right\|^{2} +\frac{3d_yL_y^2\mu_2^2}{4}.
				\label{G_y2_block}
			\end{align}
			By combing \eqref{xGk_block_E} and \eqref{G_y2_block} and rearranging terms, we obtain
			\begin{align}\label{criterion_b}
				\left\|\tilde{\mathbf{G}}(x_{t},y_{t})\right\|^{2}
				\leq\left(\frac{2}{\alpha_t^2}+3L_{y}^{2} \right) \left\|x_{t+1}-x_{t}\right\|^{2} +\frac{3}{\beta^{2}}\left\|y_{t+1}-y_{t}\right\|^{2}+ \frac{d_xL_x^2\mu_{1,t}^2}{2} +\frac{3d_yL_y^2\mu_2^2}{4}.
			\end{align}
			Denote $\theta_{t}=\frac{1}{\max\{\theta m_{t},\frac{54}{7\beta}\}}$. It is easy to calculate that
			$\frac{\frac{2}{\alpha_t^2}+3L_{y}^{2}}{m_{t}^2} \leq \theta.$
			By multiplying both sides of \eqref{criterion_b} by $\theta_{t}$ and denoting $	\mathcal{L}_{t}=\mathcal{L}(x_{t},y_{t})$, we then obtain
			\begin{align}\label{xi_t_b}
				&\theta_{t}\|\tilde{\mathbf{G}}(x_{t},y_{t})\|^{2}\notag\\
				\leq& \theta_{t}\left(\frac{2}{\alpha_t^2}+3L_{y}^{2} \right) \left\|x_{t+1}-x_{t}\right\|^{2} +\theta_{t}\frac{3}{\beta^{2}}\left\|y_{t+1}-y_{t}\right\|^{2}+\theta_{t}\left(\frac{d_xL_x^2\mu_{1,t}^2}{2} +\frac{3d_yL_y^2\mu_2^2}{4}\right)\notag\\
				\leq& m_{t}\left\|x_{t+1}-x_{t}\right\|^{2} + \frac{7}{18\beta}\left\|y_{t+1}-y_{t}\right\|^{2}
				+\theta_{t}\left(\frac{d_xL_x^2\mu_{1,t}^2}{2} +\frac{3d_yL_y^2\mu_2^2}{4}\right)\notag\\
				\leq& 	\mathcal{L}_{t}-	\mathcal{L}_{t+1} +\frac{\lambda_{t-1}-\lambda_{t}}{2}\sigma_{y}^{2}
				+\frac{16}{\beta}\left(\frac{\lambda_{t-2}}{\lambda_{t-1}}-\frac{\lambda_{t-1}}{\lambda_{t}}\right)
				\sigma_{y}^{2} 
				+r_{t}d_{y}L_{y}^{2}\mu_{2}^{2}
				\nonumber\\
				&+\theta_{t}\left(\frac{d_xL_x^2\mu_{1,t}^2}{2} +\frac{3d_yL_y^2\mu_2^2}{4}\right)+\frac{\beta d_xL_x^2\mu_{1,t}^2}{8},
			\end{align}
			where the last inequality holds by \eqref{potient_block}.
			Denote $\tilde{T}=\min\{t \mid
			\|\tilde{\mathbf{G}}(x_{t},y_{t})\|^{2}\leq\frac{\varepsilon^{2}}{4}\}$. By summarizing both sides of \eqref{xi_t_b} from $t=3$ to $t=\tilde{T}$, we then obtain
			\begin{align}\label{finasum_b}
				&\sum_{t=3}^{\tilde{T}}\theta_{t}
				\|\tilde{\mathbf{G}}(x_{t},y_{t})\|^{2}\notag\\
				\leq& 	\mathcal{L}_{3}-	\mathcal{L}_{\tilde{T}+1}
				+\frac{\lambda_{2}-\lambda_{\tilde{T}}}{2}\sigma_{y}^{2}
				+\frac{16}{\beta}\left(\frac{\lambda_{1}}{\lambda_{2}}-\frac{\lambda_{\tilde{T}-1}}{\lambda_{\tilde{T}}}\right)
				\sigma_{y}^{2}
				+\sum_{t=3}^{\tilde{T}}r_{t}d_{y}L_{y}^{2}\mu_{2}^{2} 
				\notag\\
				&+\sum_{t=3}^{\tilde{T}}\theta_{t}\left(\frac{d_xL_x^2\mu_{1,t}^2}{2} +\frac{3d_yL_y^2\mu_2^2}{4}\right)+\sum_{t=3}^{\tilde{T}}\frac{\beta d_xL_x^2\mu_{1,t}^2}{8}\nonumber\\
				\leq& 	\mathcal{L}_{3}-	\mathcal{L}_{\tilde{T}+1}
				+\frac{\lambda_{2}}{2}\sigma_{y}^{2}
				+\frac{16\lambda_{1}}{\beta\lambda_{2}}
				\sigma_{y}^{2}
				+\sum_{t=3}^{\tilde{T}}r_{t}d_{y}L_{y}^{2}\mu_{2}^{2} 
				+\sum_{t=3}^{\tilde{T}}\theta_{t}\left(\frac{d_xL_x^2\mu_{1,t}^2}{2} +\frac{3d_yL_y^2\mu_2^2}{4}\right)\nonumber\\
				&+\sum_{t=3}^{\tilde{T}}\frac{\beta d_xL_x^2\mu_{1,t}^2}{8} \notag\\
				\leq& C +\sum_{t=3}^{\tilde{T}}r_{t}d_{y}L_{y}^{2}\mu_{2}^{2}+\sum_{t=3}^{\tilde{T}}\frac{\beta d_xL_x^2\mu_{1,t}^2}{8} 
				+\sum_{t=3}^{\tilde{T}}\theta_{t}\left(\frac{d_xL_x^2\mu_{1,t}^2}{2} +\frac{3d_yL_y^2\mu_2^2}{4}\right),
			\end{align}
			where the last inequality is because $	\mathcal{L}_{\tilde{T}+1}\geq \underline{f}-\frac{56}{\beta}\sigma_{y}^{2}$ and $	\mathcal{L}_{3}\leq \bar{f}+\frac{160}{3\beta}\sigma_{y}^{2}$, and then $	\mathcal{L}_{3}-	\mathcal{L}_{\tilde{T}+1} + \left( \frac{9}{40\times 2^{1/4}\beta}+\frac{16}{\beta}\times2^{1/4} \right)\sigma_{y}^{2}\leq\bar{f}-\underline{f}+ \frac{149}{\beta}\sigma_{y}^{2} = C$. Denote  $s=\max\{\theta, \frac{54}{7\beta m_{1}}\}$. Since $m_{t}$ is an increasing sequence,  $s\geq\max\{\theta, \frac{54}{7\beta m_{t}}\}= \frac{1}{\theta_{t} m_{t}}$, which implies that $\theta_{t}\geq \frac{1}{s m_{t}}$.  It then can be calculated that
			\begin{align}\label{theta_alpha_bound_b}
				\frac{C_{1}}{s t^{1/2}}\leq\theta_{t}\le C_2.
			\end{align}
			By the definition of $\lambda_{t}$, $\beta\leq \frac{1}{10L_{y}}$ and $t\geq 1$, we also give upper bound for $r_t$ which are both defined as in \eqref{p_r_b} as follows,
			\begin{align}\label{r_bound_b}
				r_{t}\leq \frac{1}{4L_{y}}+\frac{64}{9L_{y}}t^{1/4} +\frac{10240}{81L_{y}}t^{1/2} \leq \frac{10837}{81L_{y}}t^{1/2}.
			\end{align}
			By using the definition of $\tilde{T}$ and plugging $\mu_{1,t}$ and $\mu_{2}$ that defined in \eqref{mu2_t_b} respectively, and \eqref{theta_alpha_bound_b}, \eqref{r_bound_b} into \eqref{finasum_b}, we conclude that
			\begin{align*}
				\frac{\varepsilon^{2}}{4}
				&\leq \frac{C}{\sum_{t=3}^{\tilde{T}}\theta_{t}} +\frac{\sum_{t=3}^{\tilde{T}}r_{t}}{\sum_{t=3}^{\tilde{T}}\theta_{t}}d_{y}L_{y}^{2}\mu_{2}^{2} 
				+\frac{\sum_{t=3}^{\tilde{T}}\theta_t\mu_{1,t}^2}{2\sum_{t=3}^{\tilde{T}}\theta_{t}}d_xL_x^2+\frac{\sum_{t=3}^{\tilde{T}}\mu_{1,t}^2}{8\sum_{t=3}^{\tilde{T}}\theta_{t}}\beta d_xL_x^2 +\frac{3d_yL_y^2\mu_2^2}{4}\nonumber\\
				&\le \frac{Cs}{C_1(\sqrt{\tilde{T}}-2)}+\frac{10837s}{81L_yC_1}d_{y}L_{y}^{2}\mu_{2}^{2} +\frac{3d_yL_y^2\mu_2^2}{4}+\frac{\sum_{t=3}^{\tilde{T}}\mu_{1,t}^2}{8\sum_{t=3}^{\tilde{T}}\theta_{t}}\left(4C_2+\beta\right) d_xL_x^2 \nonumber\\
				&= \frac{Cs}{C_1(\sqrt{\tilde{T}}-2)}+\frac{\varepsilon^2}{8},
			\end{align*}
			or equivalently, $\tilde{T}\leq \left(\frac{8sC}{\varepsilon^{2}C_{1}}+2\right)^{2}$.
			If $t\geq \frac{ 9^{4}\sigma_{y}^{4}}{10^{4}\beta^{4}\varepsilon^{4}}$, then $\lambda_{t}^{2}\|y_{t}\|^{2}\leq \frac{\varepsilon^{2}}{4}$. It is easy to know that $\|\mathbf{G}(x_{t},y_{t})\|^{2} \leq 2\|\tilde{\mathbf{G}}(x_{t},y_{t})\|^{2}
			+2\lambda_{t}^{2}\|y_{t}\|^{2}$. Therefore, there exists a
			$t \leq \max\left\{\left(\frac{8sC}{\varepsilon^{2}C_{1}}+2\right)^{2}, \frac{9^{4}\sigma_{y}^{4}}{10^{4}\beta^{4}\varepsilon^{4}}\right\}$
			such that $ \|\mathbf{G}(x_{t},y_{t})\|\leq \varepsilon.$	
		\end{proof}

		\begin{remark}\label{remark:1}
			If we set $\beta = \frac{1}{10L_{y}}$,
			From Theorem \ref{thm:main_b},  by plugging $\beta$ and other constants, we can conclude that for any given $\varepsilon \in (0,1)$, $T\leq \max\{(\frac{8sC}{\varepsilon^{2}C_{1}}+2)^{2}, \frac{ 9^{4}\sigma_{y}^{4}L_{y}^{4}}{\varepsilon^{4}}+1\}= \mathcal{O}\left( \varepsilon^{-4} \right)$. 
			This implies that the total complexity of the proposed ZO-AGP algorithm to obtain an $\varepsilon$-stationary point for nonconvex-concave minimax problems is bounded by $\mathcal{O}\left( (d_x+d_y)\varepsilon^{-4} \right)$.
		\end{remark}

		\subsection{A Zeroth-order Block Alternating Proximal Gradient Algorithm}\label{section Analysis_block}
		We consider a more general block-wise nonsmooth nonconvex-concave minimax problem as follows.
		\begin{equation}\tag{BP}\label{problem_block}
			\min_{x \in \mathcal{X}} \max_{y \in \mathcal{Y}}l\left( x,y \right):= f\left( x^{1},x^{2},\cdots,x^{K},y \right)+\sum_{k=1}^{K}h_{k}\left( x^{k} \right) - g\left( y \right),
		\end{equation}
		\begin{equation*}
			s.t.\quad x^{k} \in \mathcal{X}_{k},\quad k=1,\cdots,K,
		\end{equation*}
		where $f: \mathbb{R}^{d_{x}K+d_{y}}\rightarrow \mathbb{R}$ is a continuously differentiable function, $h_{k}:\mathbb{R}^{d_{x}}\rightarrow \mathbb{R}$ and $g: \mathbb{R}^{d_{y}}\rightarrow \mathbb{R}$ are some convex continuous possibly non-smooth functions, $x=\left[ x^{1},x^{2},\cdots,x^{K} \right]\in \mathcal{X}=\mathcal{X}_{1}\times \mathcal{X}_{2}\times \cdots \times\mathcal{X}_{K}\subset \mathbb{R}^{d_{x}K}$, $\mathcal{X}_{k}(k=1,\cdots,K)$ and $\mathcal{Y}$ are nonempty compact convex sets. For simplicity, we'd like to write $f\left( x^{1}, x^{2}, \cdots,x^{K}, y \right)$ as $f\left( x,y \right)$, when it's unnecessary to consider block-side below.  This composite optimization problem \eqref{problem_block} with non-smooth terms and multi-blocks structure has a wide range of applications in machine learning, signal processing and communication, 
 including distributed nonconvex optimization \cite{hajinezhad2019perturbed,liao2015semi,mateos2010distributed}, robust learning over multiple domains \cite{qian2019robust}, power control and transceiver design problem \cite{foschini1993simple,li2015multicell}, power control in the presence of a jammer \cite{gohary2009generalized}, etc. Motivated by these applications, it is of interest to develop efficient  algorithms for solving these problems with theoretical convergence  guarantees.
		In this subsection, we focus on the zeroth-order optimization algorithms to solve \eqref{problem_block}   with a black-box setting.
		
		Similar to the idea of ZO-AGP algorithm, we propose a zeroth-order block alternating proximal gradient algorithm (ZO-BAPG) to solve \eqref{problem_block}  and analyze its iteration complexity. For notational simplicity, we denote 
		\begin{align*}
			w_{t+1}^{k}&=\left[x_{t+1}^{1}, \cdots, x_{t+1}^{k-1},x_{t}^{k}, \cdots, x_{t}^{K} \right] \in \mathbb{R}^{d_{x}K},\\
			w_{t+1}^{-k}&=\left[x_{t+1}^{1}, \cdots, x_{t+1}^{k-1},x_{t}^{k+1}, \cdots, x_{t}^{K} \right] \in \mathbb{R}^{d_{x}(K-1)}.
		\end{align*}
		Let $\bar{\gamma}=\max\{\gamma_{k}\geq 0|k=1,\cdots,K\}$ and $\underline{\gamma}=\min\{\gamma_{k}\geq 0|k=1,\cdots,K\}$. Define
		\begin{align}
			\operatorname{Prox}_{\mathcal{X}_{k}} \left( w^{k} \right) &= \arg\min_{x^{k}\in \mathcal{X}_{k}}h_{k} \left( x^{k} \right) + \frac{\tau_{t}+\gamma_{k}}{2} \left\|x^{k}-w^{k}\right\|^{2},\label{prox_X}\\
			\operatorname{Prox}_{\mathcal{Y}} \left( z \right) &= \arg\max_{y\in \mathcal{Y}} - g \left( y \right) - \frac{1}{2\rho} \left\|y-z \right\|^{2}.\label{prox_Y}
		\end{align}
		Similar to the idea of ZO-AGP algorithm, instead of the original function $l\left(x,y \right) $, the proposed ZO-BAPG algorithm uses the gradient of a regularized version of the original function, i.e., 
		$\tilde{l}_t\left( x,y \right):= \tilde{f}_t\left(x,y \right)+\sum_{k=1}^{K}h_{k}\left( x^{k} \right) - g\left( y \right)$,
		where $\tilde{f}_t\left(x,y \right)=f\left(x,y \right) - \frac{\tilde{\lambda}_{t}}{2} \|y\|^{2}$ with regularization parameters $\tilde{\lambda}_t \ge 0$. Similarly, each iteration of the novel proposed ZO-BAPG algorithm conducts two proximal point ``gradient" steps for updating both $x$ and $y$, where the ``gradient" is replaced by gradient estimators shown in \eqref{x_esti} and \eqref{y_esti} respectively. More detailedly, it updates $x_t$ by minimizing a linearized approximation of $\tilde{l}\left( x,  y_t\right)$ with the gradient at point $\left( w_{t+1}^{k},y_t \right), k=1,\cdots,K$, i.e.,
		\begin{align*}
			x_{t+1}^{k}
			&=\arg\min_{x^{k}\in \mathcal{X}_{k}}\left\langle \widehat{\nabla}_{x^{k}} \tilde{f}_t\left( w_{t+1}^{k},y_{t} \right), x^{k}-x_{t}^{k}  \right\rangle + h_{k} \left( x^{k} \right) + \frac{\tau_{t}+\gamma_{k}}{2}  \left\|x^{k}-x_{t}^{k}\right\|^{2}\notag\\
			&=\operatorname{Prox}_{\mathcal{X}_{k}}\left( x_{t}^{k} - \frac{1}{\tau_{t}+\gamma_{k}} \widehat{\nabla}_{x^{k}} \tilde{f}_t\left( w_{t+1}^{k},y_{t} \right)\right),
		\end{align*}
		where $\operatorname{Prox}_{\mathcal{X}_{k}}$ is the proximal point operator onto $\mathcal{X}_{k}$ defined as in \eqref{prox_X} and 
		\begin{align}\label{x_esti_tb}
			\widehat{\nabla}_{x^{k}} \tilde{f}_t\left( w_{t+1}^{k},y_{t} \right)= \sum_{i=1}^{d_x} \frac{\Delta_{x^k} \tilde{f}_{t,i}}{\bar{\mu}_{1_{k},t}} u_{i},
		\end{align}
		with $\Delta_{x^k} \tilde{f}_{t,i} = f (x_{t+1}^{1},\cdots,x_{t+1}^{k-1},x_{t}^{k}+\bar{\mu}_{1_{k},t} u_{i},x_{t}^{k+1},\cdots,x_{t}^{K},y) - f(w_{t+1}^{k},y)$. 
		For any given $k$, $\{\bar{\mu}_{1_{k},t}\}$ is a decreasing sequence of smoothing parameters.
		For simplicity, we set $\bar{\mu}_{1_{k},t} = \bar{\mu}_{1,t}$ $(k=1,\cdots,K)$.
		Similarly, it updates $y_t$ by maximizing a linearized approximation of $\tilde{l}\left( x_{t+1}, y\right)$ minus some regularized term where the gradient at point $\left( x_{t+1},y_t \right)$ is replaced by a gradient estimator shown in  \eqref{y_esti}, i.e.,
		\begin{align*}
			y_{t+1}
			&= \arg\max_{y\in \mathcal{Y}}\left\langle \widehat{\nabla}_{y} \tilde{f}_t\left( x_{t+1},y_{t} \right), y-y_{t} \right\rangle - g \left( y \right) - \frac{1}{2\rho} \left\|y-y_{t} \right\|^{2}\notag\\
			&=\operatorname{Prox}_{\mathcal{Y}} \left( y_{t} + \rho \left(\widehat{\nabla}_{y} f\left( x_{t+1},y_{t} \right)- \tilde{\lambda}_{t}y_{t}\right)\right),
		\end{align*}
		where $\operatorname{Prox}_{\mathcal{Y}}$ is the proximal point operator onto $\mathcal{Y}$ defined as in \eqref{prox_Y} and 
		\begin{align}
			\widehat{\nabla}_{y} \tilde{f}_t\left( x_{t+1},y_{t} \right)=\widehat{\nabla}_{y} f\left(x,y_{t}\right)- \tilde{\lambda}_{t}y_{t} = \sum_{i=1}^{d_y} \frac{\left[f\left(x,y_{t}
				+\bar{\mu}_{2} {v}_{i})-f(x,y_{t}\right)\right]}{\bar{\mu}_{2}} {v}_{i}- \tilde{\lambda}_{t}y_{t}.
		\end{align}
		The proposed ZO-BAPG algorithm is formally stated in Algorithm \ref{algo2}.
		\begin{algorithm}
			\caption{(ZO-BAPG Algorithm)}
			\label{algo2}
			\begin{algorithmic}
				\STATE{\textbf{Step 1}:  Input: $x_{1}^{j}, y_{1}$, $\rho$, $\lambda_{1}$, $\tau_{1}$, $\gamma_{j}$  $(j=1,\cdots,K)$; Set $t=1$.}
				\STATE{\textbf{Step 2}:~Calculate $\tau_t$, and for $k=1,\cdots,K$, calculate $\gamma_{k}$ and perform the following update for $x_{t}$:  	\begin{align}\label{algolstep_x_block}
						x_{t+1}^{k}
						&=\operatorname{Prox}_{\mathcal{X}_{k}}\left( x_{t}^{k} - \frac{1}{\tau_{t}+\gamma_{k}} \widehat{\nabla}_{x^{k}} \tilde{f}_t\left( w_{t+1}^{k},y_{t} \right)\right).
				\end{align}}		
				\STATE{\textbf{Step 3}:~Calculate $\lambda_t$, and perform the following update for $y_t$:
					\begin{align}\label{algolstep_y_block}
						y_{t+1}
						&=\operatorname{Prox}_{\mathcal{Y}} \left( y_{t} + \rho \left(\widehat{\nabla}_{y} f\left( x_{t+1},y_{t} \right)- \lambda_{t}y_{t}\right)\right).
				\end{align}}
				\STATE{\textbf{Step 4}:~If converges, stop; otherwise, set $t=t+1, $ go to Step 2.}
			\end{algorithmic}
		\end{algorithm}
		
		Before we analyze the convergence of Algorithm \ref{algo2}, we define the stationarity gap as the termination criterion as follows.
		
		\begin{definition}
			At each iteration of Algorithm \ref{algo2}, the stationarity gap for problem \eqref{problem_block} is defined as:
			\begin{equation}
				\nabla\mathcal{G}\left(x_{t},y_{t}\right) =
				\begin{bmatrix}
					\left( \tau_{t}+\gamma_{1} \right)\left(x_{t}^{1} - \operatorname{Prox}_{\mathcal{X}_{1}} \left(x_{t}^{1}- \frac{1}{\tau_{t}+\gamma_{1}} \nabla_{x^{1}}f\left(x_{t},y_{t}\right)\right)\right)\\
					\vdots\\
					\left( \tau_{t}+\gamma_{K} \right) \left(x_{t}^{K} - \operatorname{Prox}_{\mathcal{X}_{K}} \left(x_{t}^{K}- \frac{1}{\tau_{t}+\gamma_{K}}\nabla_{x^{K}}f\left(x_{t},y_{t}\right)\right)\right)\\
					\frac{1}{\rho} \left(y_{t} - \operatorname{Prox}_{\mathcal{Y}} \left(y_{t}+\rho\nabla_{y}f\left(x_{t},y_{t}\right)\right)\right)
				\end{bmatrix}.
			\end{equation}
		\end{definition}
		
		\begin{definition}
			At each iteration of Algorithm \ref{algo2}, we can also define the stationarity gap as:
			\begin{equation}\label{conver_block2}
				\nabla\tilde{\mathcal{G}}\left(x_{t},y_{t}\right) =
				\begin{bmatrix}
					\left( \tau_{t}+\gamma_{1} \right) \left(x_{t}^{1} - \operatorname{Prox}_{\mathcal{X}_{1}} \left(x_{t}^{1}- \frac{1}{\tau_{t}+\gamma_{1}} \nabla_{x^{1}}\tilde{f}_t\left(x_{t},y_{t}\right)\right)\right)\\
					\vdots\\
					\left( \tau_{t}+\gamma_{K} \right) \left(x_{t}^{K} - \operatorname{Prox}_{\mathcal{X}_{K}} \left(x_{t}^{K}- \frac{1}{\tau_{t}+\gamma_{K}} \nabla_{x^{K}}\tilde{f}_t\left(x_{t},y_{t}\right)\right)\right)\\
					\frac{1}{\rho} \left(y_{t} - \operatorname{Prox}_{\mathcal{Y}} \left(y_{t}+\rho\nabla_{y}\tilde{f}_t\left(x_{t},y_{t}\right)\right)\right)
				\end{bmatrix}.
			\end{equation}
			For simplicity, we denote $\nabla\tilde{\mathcal{G}}_{t}=\nabla\tilde{\mathbf{G}}\left(x_{t},y_{t}\right)$ and for $k=1,\cdots, K$,
			\begin{align*}
				\left(\nabla\tilde{\mathcal{G}}_{t}\right)_{x^{k}}&=\left( \tau_{t}+\gamma_{k} \right) \left(x_{t}^{k} - \operatorname{Prox}_{\mathcal{X}_{k}} \left(x_{t}^{k}-\frac{1}{\tau_{t}+\gamma_{k}} \nabla_{x^{k}}\tilde{f}_t\left(x_{t},y_{t}\right)\right)\right),\\
				\left(\nabla\tilde{\mathcal{G}}_{t}\right)_{y}&=\frac{1}{\rho}\left(y_{t}
				-\operatorname{Prox}_{\mathcal{Y}}
				\left(y_{t}+\rho\nabla_{y}\tilde{f}_t\left(x_{t},y_{t}\right)\right)\right).
			\end{align*}
		\end{definition}
		
		Firstly, we need to make the following assumption about the smoothness of $f$.
		
		\begin{assumption}\label{a4_block}
			$f\left(x,y\right)$ has Lipschitz continuous gradients, i.e., there exist $L_{x_{k}}$, $L_{y}$ such that for $k=1,\cdots,K$,
			\begin{align*}
				\left\|\nabla_{x^{k}}f\left(x_{1},y\right) - \nabla_{x^{k}}f\left(x_{2},y\right)\right\| & \leq L_{x_{k}} \left\|x_{1}-x_{2}\right\|, \quad \forall x_{1},x_{2}\in\mathcal{X}, \forall y\in \mathcal{Y}, \\
				\left\|\nabla_{y}f\left(x_{1},y\right) - \nabla_{y}f\left(x_{2},y\right)\right\| &\leq L_{y} \left\|x_{1}-x_{2}\right\|, \quad \forall x_{1},x_{2}\in\mathcal{X},\forall y\in \mathcal{Y}, \\
				\left\|\nabla_{y}f\left(x,y_{1}\right) - \nabla_{y}f\left(x,y_{2}\right)\right\| &\leq L_{y} \left\|y_{1}-y_{2}\right\|, \quad \forall x\in \mathcal{X}, \forall y_{1},y_{2}\in\mathcal{Y}.
			\end{align*}
			Denote $L_{x}=\max\left\{ L_{x_{1}},\cdots,L_{x_{K}} \right\}$.
		\end{assumption}
		
		Since that $f\left(x,y\right)$, $h_{k}:\mathbb{R}^{d_{x}}\rightarrow \mathbb{R}$$(k=1,\cdots,K)$ and $g: \mathbb{R}^{d_{y}}\rightarrow \mathbb{R}$ are some  continuous functions, and $\mathcal{X}_{k}$ $(k=1,\cdots,K),$ $\mathcal{Y}$ are  nonempty compact convex sets. Then, $l(x,y)$ is bounded, i.e., there exist finite numbers $\bar{l}$ and $\underline{l}$, such that $\underline{l}\leq l(x,y) \leq \bar{l}$. Now, we are ready to show the descent property as follows.
		
		\begin{lemma}\label{dec_y}
			Suppose that Assumption \ref{a4_block} and $\rho \leq \frac{1}{L_{y}+2\lambda_{1}}$ hold. Assume the sequence $\{(x_{t},y_{t})\}$ is generated by Algorithm \ref{algo2}. Define $	\mathcal{L}(x_{t},y_{t})=l(x_{t},y_{t})+\frac{16}{\rho^{2}\tilde{\lambda}_{t}}\|y_{t}-y_{t-1}\|^{2}
			-\frac{16}{\rho}\left(\frac{\tilde{\lambda}_{t-2}}{\tilde{\lambda}_{t-1}}-1\right)\|y_{t}\|^{2} -\frac{\tilde{\lambda}_{t-1}}{2}\|y_{t}\|^{2}
			-\frac{11}{\rho}\|y_{t}-y_{t-1}\|^{2}.$
			For any $t\geq 3$, we have
			\begin{align*}
				\mathcal{L}(x_{t+1},y_{t+1})\leq 	& \mathcal{L} (x_{t},y_{t})
				-\tilde{c}_{t}\|y_{t+1}-y_{t}\|^{2}
				-\tilde{m}_{t}\|x_{t+1}-x_{t}\|^{2}+\frac{K\rho d_xL_x^2 \bar{\mu}_{1,t}^2}{8}\\
				&+ \frac{\tilde{\lambda}_{t-1}-\tilde{\lambda}_{t}}{2}\|y_{t+1}\|^{2}
				+\frac{16}{\rho}\left(\frac{\tilde{\lambda}_{t-2}}{\tilde{\lambda}_{t-1}}-\frac{\tilde{\lambda}_{t-1}}{\tilde{\lambda}_{t}}\right)
				\|y_{t}\|^{2}
				+ \tilde{r}_{t}d_{y}L_{y}^{2}\bar{\mu}_{2}^{2},
			\end{align*}
			where the parameters are chosen as follows,
			\begin{align}
				\tilde{c}_{t} &= \frac{15}{2\rho}-\frac{16}{\rho^{2}}\left(\frac{1}{\tilde{\lambda}_{t+1}}
				-\frac{1}{\tilde{\lambda}_{t}}\right),\quad \tilde{m}_{t} = \tau_{t} + \underline{\gamma}-\left(\frac{L_{x}}{2}+\frac{1}{2\rho}+\frac{384}{\rho\tilde{\lambda}_{t}^{2}}L_{y}^{2}+\frac{3\rho L_{y}^{2}}{2}\right),\label{nc_m_b}\\
				\tilde{r}_{t} &= \frac{19\rho}{8}+\frac{32}{\tilde{\lambda}_{t}}+\frac{256}{\rho\tilde{\lambda}_{t}^{2}}.\label{np_r_b}
			\end{align}
			
		\end{lemma}
		
		\begin{proof}
			By Assumption \ref{a4_block} and the convexity of $h_{k}$, $f(x,y)$ has Lipschtiz continuous gradient w.r.t. $x$ , which implies that
			\begin{align}\label{l6(0)}
				&l(x_{t+1}^{k},w_{t+1}^{-k},y_{t})- l(w_{t+1}^{k},y_{t})\notag\\
				\leq& \left\langle \nabla_{x^{k}} f(w_{t+1}^{k},y_{t}) + \partial h_{k} \left( x_{t+1}^{k} \right), x_{t+1}^{k}-x_{t}^{k} \right\rangle	+ \frac{L_{x}}{2}\|x_{t+1}^{k}-x_{t}^{k}\|^{2}\notag\\
				=& \left\langle \widehat{\nabla}_{x^{k}} \tilde{f}_t(w_{t+1}^{k},y_{t}) + \partial h_{k} \left( x_{t+1}^{k} \right), x_{t+1}^{k}-x_{t}^{k} \right\rangle	+ \frac{L_{x}}{2}\|x_{t+1}^{k}-x_{t}^{k}\|^{2}\notag\\
				&+ \left\langle \nabla_{x^{k}} f(w_{t+1}^{k},y_{t}) -  \widehat{\nabla}_{x^{k}} \tilde{f}_t(w_{t+1}^{k},y_{t}) , x_{t+1}^{k}-x_{t}^{k} \right\rangle.
			\end{align}
			Firstly, we estimate the second term in the r.h.s of \eqref{l6(0)}. From the optimality of $x_{t+1}$'s update in \eqref{algolstep_x_block} of  Algorithm \ref{algo2},  $\forall x^{k}\in\mathcal{X}_{k}$, we obtain
			\begin{equation}
				\left\langle \widehat{\nabla}_{x^{k}} \tilde{f}_t(w_{t+1}^{k},y_{t}) + \partial h_{k} ( x_{t+1}^{k} ) + (\tau_{t} + \gamma_{k})(x_{t+1}^{k} - x_{t}^{k}), x^{k}-x_{t+1}^{k} \right\rangle \geq 0.\label{l6_2nd1}
			\end{equation}
			By setting $x^{k}=x_{t}^{k}$ and rearranging the terms, we have
			\begin{equation}
				\left\langle \widehat{\nabla}_{x^{k}} \tilde{f}_t(w_{t+1}^{k},y_{t}) + \partial h_{k} \left( x_{t+1}^{k} \right),x_{t+1}^{k}-x_{t}^{k} \right\rangle \leq - (\tau_{t} + \gamma_{k})\left\| x_{t+1}^{k}-x_{t}^{k} \right\|^{2}.\label{l6_2nd2}
			\end{equation}
			On the other hand, we estimate the last term in the right side of \eqref{l6(0)}. By the Cauchy-Schwarz inequality and Lemma \ref{varboundlemma_b}, we have
			\begin{align}
				\left\langle \nabla_{x^{k}} f(w_{t+1}^{k},y_{t}) -  \widehat{\nabla}_{x^{k}} \tilde{f}_t(w_{t+1}^{k},y_{t}) , x_{t+1}^{k}-x_{t}^{k} \right\rangle
				\le\frac{1}{2\rho}\|x_{t+1}^k-x_{t}^k\|^{2}+\frac{\rho d_xL_x^2 \bar{\mu}_{1,t}^2}{8}.\label{l6_4th2_bl}
			\end{align}
			By plugging \eqref{l6_2nd2} and \eqref{l6_4th2_bl} into \eqref{l6(0)} and summing up it from $k=1$ to $K$, then by the definition of $\underline{\gamma}$, we obtian
			\begin{equation}
				l(x_{t+1},y_{t})\leq l(x_{t},y_{t})- \left(\tau_{t} + \underline{\gamma}- \frac{L_x}{2}-\frac{1}{2\rho}\right)\left\|x_{t+1} - x_{t}\right\|^{2} +\frac{K\rho d_xL_x^2 \bar{\mu}_{1,t}^2}{8}.\label{decentx2}
			\end{equation}
			The rest of the proof is almost the same as Lemma \ref{dislemma_block} and \ref{decentxy_block} since the ascent step of $y$'s update can be similarly estimated.
			By replacing $\nabla _yf(x_{k+1},y_k)$ with $\nabla _yf(x_{k+1},y_k)-\partial g(y_{k+1})$, and using $\langle \partial g(y_{k+1})-\partial g(y_k),y_{k+1}-y_k \rangle \geq0$. The proof is completed  by some parameters replacement, e.g., $\beta$ with $\rho$, $\lambda$ with $\tilde{\lambda}$, $c_t$ with $\tilde{c}_t$, $m_t$ with $\tilde{m}_t$, $r_t$ with $\tilde{r}_t$ and $\mu_2$ with $\bar{\mu}_2$.
		\end{proof}

		\begin{theorem}\label{thm:main_block}
			Suppose that Assumption \ref{a4_block} holds. The sequence $\{(x_{t},y_{t})\}$ is generated by Algorithm \ref{algo2}. Let $\tilde{\lambda}_{t}=\frac{9}{20\rho t^{1/4}}, \quad \rho\leq\frac{1}{10L_{y}},\quad \tau_{t} =\frac{L_{x}}{2}+\frac{1}{2\rho}+\frac{768}{\rho\tilde{\lambda}_{t}^{2}}L_{y}^{2}+\frac{3\rho L_{y}^{2}}{2}- \underline{\gamma}.$
			Let $\theta =24+\frac{9^{4}\left[6(\frac{L_{x}}{2}+\frac{1}{2\rho}+\frac{3\rho L_{y}^{2}}{2}-\underline{\gamma}+\bar{\gamma} )^{2} +3L_{y}^{2}+3KL_x^2\right]}{\left(384\times 20^{2}\rho L_{y}^{2}\right)^{2}}$,
			$s =\max\left\{\theta, \frac{27^{2}}{448\times 20^{2}\rho^{2}L_{y}^{2}}\right\}$, $C_{1}=\frac{27}{128\times20^{2}\rho L_{y}^{2}}$, $C_2=\max\{\frac{7\rho}{54}, C_1\}$.
			Then for any $t\geq 3$ and for any given $\varepsilon \in (0,1)$, if we set
			\begin{align}
				\bar{\mu}_{1,t}=\frac{\varepsilon}{L_xt^{1/4}}\sqrt{\frac{C_1}{2s(6C_2+\rho)Kd_x}},\quad
				\bar{\mu}_{2}=\frac{\varepsilon}{4L_{y}\sqrt{d_y}}\sqrt{\frac{1}{\frac{10837s}{81L_yC_1}+\frac{3}{4}}},\label{mu2_t_block}
			\end{align}
			then
			$$T(\varepsilon) \leq \max\left\{\left(\frac{8sC}{\varepsilon^{2}C_{1}}+2\right)^{2}, \frac{ 9^{4}\sigma_{y}^{4}}{10^{4}\rho^{4}\varepsilon^{4}}\right\},$$
			where $C=\bar{l}-\underline{l}+ \frac{149}{\rho}\sigma_{y}^{2} $ and $\sigma_{y}=\max\{\|y\||y\in\mathcal{Y}\}$.
		\end{theorem}
		
		\begin{proof}
			Form the setting of $\tilde{\lambda}_{t}$ and $\rho$, it holds that
			$\frac{1}{\tilde{\lambda}_{t+1}}-\frac{1}{\tilde{\lambda}_{t}}\leq\frac{4\rho}{9}$ and  $\tau_{t}+\underline{\gamma}>\frac{L_{x}}{2}+\frac{1}{2\rho}+\frac{384}{\rho\tilde{\lambda}_{t}^{2}}L_{y}^{2}+\frac{3\rho L_{y}^{2}}{2},$
			and $\tilde{m}_{t}=\frac{384}{\rho\tilde{\lambda}_{t}^{2}}L_{y}^{2}$. Then we have $\tilde{c}_{t}\geq\frac{7}{18\rho}$ which is defined as in \eqref{nc_m_b}.  Note that $\tilde{r}_{t}>0$ defined as in \eqref{np_r_b}. On the one hand, by Lemma \ref{dec_y}, we obtain
			\begin{align}
				\mathcal{L}(x_{t+1},y_{t+1})\leq& \mathcal{L}(x_{t},y_{t})
				-\frac{7}{18\rho}\|y_{t+1}-y_{t}\|^{2}
				-\tilde{m}_{t}\|x_{t+1}-x_{t}\|^{2}+\frac{K\rho d_xL_x^2\bar{\mu}_{1,t}^2}{8}\notag\\
				&
				+\frac{\tilde{\lambda}_{t-1}-\tilde{\lambda}_{t}}{2}\sigma_{y}^{2}
				+\frac{16}{\rho}\left(\frac{\tilde{\lambda}_{t-2}}{\tilde{\lambda}_{t-1}}-\frac{\tilde{\lambda}_{t-1}}{\tilde{\lambda}_{t}}\right)
				\sigma_{y}^{2}+\tilde{r}_{t}d_{y}L_{y}^{2}\bar{\mu}_{2}^{2},\label{potient_b}
			\end{align}
			where $\frac{\tilde{\lambda}_{t-1}-\tilde{\lambda}_{t}}{2} >0$ and $\frac{16}{\rho}\left(\frac{\tilde{\lambda}_{t-2}}{\tilde{\lambda}_{t-1}}-\frac{\tilde{\lambda}_{t-1}}{\tilde{\lambda}_{t}}\right) >0$.
			On the other hand, by nonexpansiveness of the proximal operator and using the Cauchy-Schwarz inequality, we have
			\begin{align}\label{xGk_b}
				\left\|(\nabla\tilde{\mathcal{G}}_{t})_{x^{k}}\right\|
				&=(\tau_{t}+\gamma_{k}) \left\|x_{t+1}^{k}-x_{t}^{k} - (x_{t+1}^{k}-\operatorname{Prox}_{\mathcal{X}_{k}}(x_{t}^{k}-\frac{1}{\tau_{t}+\gamma_{k}}\nabla_{x^{k}}
				f(x_{t},y_{t})))\right\|\notag\\
				&\leq(\tau_{t}+\gamma_{k})\left\|x_{t+1}^{k}-x_{t}^{k}\right\|
				+\left\|\widehat{\nabla}_{x^{k}} \tilde{f}_t(w_{t+1}^{k},y_{t})
				-\nabla_{x^{k}}f(x_{t},y_{t})\right\|\nonumber\\
				&\leq(\tau_{t}+\gamma_{k})\left\|x_{t+1}^{k}-x_{t}^{k}\right\|
				+\left\|\widehat{\nabla}_{x^{k}} \tilde{f}_t(w_{t+1}^{k},y_{t})
				-\nabla_{x^{k}}f(w_{t+1}^{k},y_{t})\right\|\nonumber\\
				&\quad+\left\|\nabla_{x^{k}}f(w_{t+1}^{k},y_{t})-\nabla_{x^{k}}f(x_{t},y_{t})\right\|,
			\end{align}
			and
			\begin{align}\label{yG_block}
				\left\|(\nabla\tilde{\mathcal{G}}_{t})_{y}\right\|
				&\leq\frac{1}{\rho}\left\|y_{t+1}-y_{t}\right\|
				+\frac{1}{\rho}\left\|y_{t+1}-\operatorname{Prox}_{\mathcal{Y}}(y_{t}+\rho\nabla_{y}
				\tilde{f}_t(x_{t},y_{t}))\right\|\notag\\
				&\leq\frac{1}{\rho}\left\|y_{t+1}-y_{t}\right\|
				+\left\|\widehat{\nabla}_{y}f(x_{t+1},y_{t})-\tilde{\lambda}_{t}y_{t}
				-(\nabla_{y}f(x_{t},y_{t})-\tilde{\lambda}_{t}y_{t})\right\|\notag\\
				&\leq\frac{1}{\rho}\left\|y_{t+1}-y_{t}\right\|
				+\left\|\widehat{\nabla}_{y}f(x_{t+1},y_{t})-\nabla_{y}
				f(x_{t},y_{t})\right\|\nonumber\\
				&\le \frac{1}{\rho}\left\|y_{t+1}-y_{t}\right\|
				+\left\|\widehat{\nabla}_{y}f(x_{t+1},y_{t})-\nabla_{y}
				f(x_{t+1},y_{t})\right\|\nonumber\\
				&\quad+\left\|\nabla_{y}f(x_{t+1},y_{t})-\nabla_{y}
				f(x_{t},y_{t})\right\|.
			\end{align}    
			Next, by using the Cauchy-Schwarz inequality, Lemma \ref{varboundlemma_b}, the definition of $\bar{\gamma}$ and the fact that $\left\|w_{t+1}^{k}-x_{t}\right\| \leq \left\|x_{t+1}-x_{t}\right\|$, we obtain
			\begin{align}\label{xGk_b_E}
				\left\|(\nabla\tilde{\mathcal{G}}_{t})_{x^k}\right\|^{2}
				\leq& 3(\tau_{t}+\gamma_{k})^2\left\|x_{t+1}^k-x_{t}^k\right\|^{2} + \frac{3d_xL_x^2\bar{\mu}_{1,t}^2}{4}+3L_x^2\|x_{t+1}-x_t\|^2 \nonumber\\		
				\leq& 3(\tau_{t}+\bar{\gamma})^2\left\|x_{t+1}^k-x_{t}^k\right\|^{2} + \frac{3d_xL_x^2\bar{\mu}_{1,t}^2}{4}+3L_x^2\|x_{t+1}-x_t\|^2.
			\end{align}
			Summing up \eqref{xGk_b_E} from $k=1$ to $K$, then we have
			\begin{align}\label{xGk_b_En}
				\left\|(\nabla\tilde{\mathcal{G}}_{t})_{x}\right\|^{2}
				\leq 3\left((\tau_{t}+\bar{\gamma})^2+KL_x^2\right)\left\|x_{t+1}-x_{t}\right\|^{2} + \frac{3Kd_xL_x^2\bar{\mu}_{1,t}^2}{4} .
			\end{align}
			Similarly, we have
			\begin{align}
				\left\|(\nabla\tilde{\mathcal{G}}_{t})_{y}\right\|^{2}
				\leq \frac{3}{\rho^{2}}\left\|y_{t+1}-y_{t}\right\|^{2}+3L_y^2\left\|x_{t+1}-x_{t}\right\|^{2} +\frac{3d_yL_y^2\bar{\mu}_2^2}{4}.
				\label{G_y2_b}
			\end{align}
			By combing \eqref{xGk_b_En} and \eqref{G_y2_b} and rearranging terms, we obtain
			\begin{align}
				\left\|\tilde{\mathcal{G}}(x_{t},y_{t})\right\|^{2}
				\leq&\left(3(\tau_{t}+\bar{\gamma})^2+3KL_x^2+3L_{y}^{2} \right) \left\|x_{t+1}-x_{t}\right\|^{2} +\frac{3}{\rho^{2}}\left\|y_{t+1}-y_{t}\right\|^{2}\nonumber\\ &+ \frac{3Kd_xL_x^2\bar{\mu}_{1,t}^2}{4}+\frac{3d_yL_y^2\bar{\mu}_2^2}{4}.\label{criterion_block}
			\end{align}
			Denote $\theta_{t}=\frac{1}{\max\{\theta \tilde{m}_{t},\frac{54}{7\rho}\}}$. It is easy to calculate that
			$\frac{3(\tau_{t}+\bar{\gamma})^2+3KL_x^2+3L_{y}^{2}}{\tilde{m}_{t}^2} \leq \theta.$
			By multiplying both sides of \eqref{criterion_block} by $\theta_{t}$ and denoting $\mathcal{L}_{t}=\mathcal{L}(x_{t},y_{t})$, we then obtain
			\begin{align}\label{xi_t_block}
				&\theta_{t}\|\tilde{\mathcal{G}}(x_{t},y_{t})\|^{2}\notag\\
				\leq& \tilde{m}_{t}\left\|x_{t+1}-x_{t}\right\|^{2} + \frac{7}{18\rho}\left\|y_{t+1}-y_{t}\right\|^{2}
				+\theta_{t}\left(\frac{3Kd_xL_x^2\bar{\mu}_{1,t}^2}{4} +\frac{3d_yL_y^2\bar{\mu}_2^2}{4}\right)\notag\\
				\leq& \mathcal{L}_{t}-\mathcal{L}_{t+1} +\frac{\tilde{\lambda}_{t-1}-\tilde{\lambda}_{t}}{2}\sigma_{y}^{2}
				+\frac{16}{\rho}\left(\frac{\tilde{\lambda}_{t-2}}{\tilde{\lambda}_{t-1}}-\frac{\tilde{\lambda}_{t-1}}{\tilde{\lambda}_{t}}\right)
				\sigma_{y}^{2} 
				+\tilde{r}_{t}d_{y}L_{y}^{2}\bar{\mu}_{2}^{2}
				\nonumber\\
				&+\theta_{t}\left(\frac{3Kd_xL_x^2\bar{\mu}_{1,t}^2}{4} +\frac{3d_yL_y^2\bar{\mu}_2^2}{4}\right)+\frac{K\rho d_xL_x^2\bar{\mu}_{1,t}^2}{8},
			\end{align}
			where the last inequality holds by \eqref{potient_b}.
			Denote $\tilde{T}=\min\{t \mid
			\|\tilde{\mathcal{G}}(x_{t},y_{t})\|^{2}\leq\frac{\varepsilon^{2}}{4}\}$. By summarizing both sides of \eqref{xi_t_block} from $t=3$ to $t=\tilde{T}$, we then obtain
			\begin{align}\label{finasum_block}
				&\sum_{t=3}^{\tilde{T}}\theta_{t}
				\|\tilde{\mathcal{G}}(x_{t},y_{t})\|^{2}\notag\\
				\leq& \mathcal{L}_{3}-\mathcal{L}_{\tilde{T}+1}
				+\frac{\tilde{\lambda}_{2}-\tilde{\lambda}_{\tilde{T}}}{2}\sigma_{y}^{2}
				+\frac{16}{\rho}\left(\frac{\tilde{\lambda}_{1}}{\tilde{\lambda}_{2}}-\frac{\tilde{\lambda}_{\tilde{T}-1}}{\tilde{\lambda}_{\tilde{T}}}\right)
				\sigma_{y}^{2}
				+\sum_{t=3}^{\tilde{T}}\tilde{r}_{t}d_{y}L_{y}^{2}\bar{\mu}_{2}^{2} 
				\notag\\
				&+\sum_{t=3}^{\tilde{T}}\theta_{t}\left(\frac{3Kd_xL_x^2\bar{\mu}_{1,t}^2}{4} +\frac{3d_yL_y^2\bar{\mu}_2^2}{4}\right) +\sum_{t=3}^{\tilde{T}}\frac{K\rho d_xL_x^2\bar{\mu}_{1,t}^2}{8}\nonumber\\
				\leq& \mathcal{L}_{3}-\mathcal{L}_{\tilde{T}+1}
				+\frac{\tilde{\lambda}_{2}}{2}\sigma_{y}^{2}
				+\frac{16\tilde{\lambda}_{1}}{\rho\tilde{\lambda}_{2}}
				\sigma_{y}^{2}
				+\sum_{t=3}^{\tilde{T}}\tilde{r}_{t}d_{y}L_{y}^{2}\bar{\mu}_{2}^{2} +\sum_{t=3}^{\tilde{T}}\frac{K\rho d_xL_x^2\bar{\mu}_{1,t}^2}{8}\nonumber\\
				&+\sum_{t=3}^{\tilde{T}}\theta_{t}\left(\frac{3Kd_xL_x^2\bar{\mu}_{1,t}^2}{4} +\frac{3d_yL_y^2\bar{\mu}_2^2}{4}\right) \notag\\
				\leq& C +\sum_{t=3}^{\tilde{T}}\tilde{r}_{t}d_{y}L_{y}^{2}\bar{\mu}_{2}^{2}+\sum_{t=3}^{\tilde{T}}\frac{K\rho d_xL_x^2\bar{\mu}_{1,t}^2}{8} 
				+\sum_{t=3}^{\tilde{T}}\theta_{t}\left(\frac{3Kd_xL_x^2\bar{\mu}_{1,t}^2}{4} +\frac{3d_yL_y^2\bar{\mu}_2^2}{4}\right),
			\end{align}
			where the last inequality is because $\mathcal{L}_{\tilde{T}+1}\geq \underline{l}-\frac{56}{\rho}\sigma_{y}^{2}$ and $\mathcal{L}_{3}\leq \bar{l}+\frac{160}{3\rho}\sigma_{y}^{2}$, and then $\mathcal{L}_{3}-\mathcal{L}_{\tilde{T}+1} + \left( \frac{9}{40\times 2^{1/4}\rho}+\frac{16}{\rho}\times2^{1/4} \right)\sigma_{y}^{2}\leq\bar{l}-\underline{l}+ \frac{149}{\rho}\sigma_{y}^{2} = C$. Denote  $s=\max\{\theta, \frac{54}{7\rho \tilde{m}_{1}}\}$. Since $m_{t}$ is an increasing sequence,  $s\geq\max\{\theta, \frac{54}{7\rho \tilde{m}_{t}}\}= \frac{1}{\theta_{t} \tilde{m}_{t}}$, which implies that $\theta_{t}\geq \frac{1}{s \tilde{m}_{t}}$.  It then can be calculated that
			\begin{align}\label{theta_alpha_bound_bl}
				\frac{C_{1}}{s t^{1/2}}\leq\theta_{t}\le C_2.
			\end{align}
			By the definition of $\tilde{\lambda}_{t}$, $\rho\leq \frac{1}{10L_{y}}$ and $t\geq 1$, we also give upper bound for $\tilde{r}_t$ which are both defined as in \eqref{np_r_b} as follows,
			\begin{align}\label{r_bound_bl}
				\tilde{r}_{t}\leq \frac{1}{4L_{y}}+\frac{64}{9L_{y}}t^{1/4} +\frac{10240}{81L_{y}}t^{1/2} \leq \frac{10837}{81L_{y}}t^{1/2}.
			\end{align}
			By using the definition of $\tilde{T}$ and plugging $\bar{\mu}_{1,t}$ and $\bar{\mu}_{2}$ that defined in \eqref{mu2_t_block} respectively, and \eqref{theta_alpha_bound_bl}, \eqref{r_bound_bl} into \eqref{finasum_block}, we conclude that
			\begin{align*}
				\frac{\varepsilon^{2}}{4}
				&\leq \frac{C}{\sum_{t=3}^{\tilde{T}}\theta_{t}} +\frac{\sum_{t=3}^{\tilde{T}}\tilde{r}_{t}}{\sum_{t=3}^{\tilde{T}}\theta_{t}}d_{y}L_{y}^{2}\bar{\mu}_{2}^{2} 
				+\frac{3\sum_{t=3}^{\tilde{T}}\theta_t\bar{\mu}_{1,t}^2}{4\sum_{t=3}^{\tilde{T}}\theta_{t}}Kd_xL_x^2+\frac{\sum_{t=3}^{\tilde{T}}\bar{\mu}_{1,t}^2}{8\sum_{t=3}^{\tilde{T}}\theta_{t}}K\rho d_xL_x^2\nonumber\\ &\quad+\frac{3d_yL_y^2\bar{\mu}_2^2}{4}\nonumber\\
				&\le \frac{Cs}{C_1(\sqrt{\tilde{T}}-2)}+\frac{10837s}{81L_yC_1}d_{y}L_{y}^{2}\bar{\mu}_{2}^{2} +\frac{3d_yL_y^2\bar{\mu}_2^2}{4}+\frac{\sum_{t=3}^{\tilde{T}}\bar{\mu}_{1,t}^2}{8\sum_{t=3}^{\tilde{T}}\theta_{t}}\left(6C_2+\rho\right) Kd_xL_x^2 \nonumber\\
				&= \frac{Cs}{C_1(\sqrt{\tilde{T}}-2)}+\frac{\varepsilon^2}{8},
			\end{align*}
			or equivalently, $\tilde{T}\leq \left(\frac{8sC}{\varepsilon^{2}C_{1}}+2\right)^{2}$.
			If $t\geq \frac{ 9^{4}\sigma_{y}^{4}}{10^{4}\rho^{4}\varepsilon^{4}}$, then $\tilde{\lambda}_{t}^{2}\|y_{t}\|^{2}\leq \frac{\varepsilon^{2}}{4}$. It is easy to know that $\|\mathcal{G}(x_{t},y_{t})\|^{2} \leq 2\|\tilde{\mathcal{G}}(x_{t},y_{t})\|^{2}
			+2\tilde{\lambda}_{t}^{2}\|y_{t}\|^{2}$. Therefore, there exists a
			$t \leq \max\left\{\left(\frac{8sC}{\varepsilon^{2}C_{1}}+2\right)^{2}, \frac{9^{4}\sigma_{y}^{4}}{10^{4}\rho^{4}\varepsilon^{4}}\right\}$
			such that $ \|\mathcal{G}(x_{t},y_{t})\|\leq \varepsilon.$	
		\end{proof}

		\begin{table}
			\centering	
			\caption{Zeroth-order algorithms for deterministic nonconvex-(strongly) concave minimax problems}
			\resizebox{\textwidth}{!}{
				\begin{threeparttable}
					\begin{tabular}{|c|c|c|c|c|}			
						\hline
						\multirow{2}{*}{Setting} & \multirow{2}{*}{Algorithm} & Iteration & Per iteration function &\tnote{1}  Total \\
						~ & ~ & complexity & value  estimations & complexity \\
						\hline
						\tnote{2} NC-SC  & ZO-Min-Max &\multirow{2}{*}{$\mathcal{O}\left( \epsilon^{-2}\right)$}  & \multirow{2}{*}{$\mathcal{O}\left(\left(d_{x}+d_{y}\right) \epsilon^{-2}\right)$} &\multirow{2}{*}{$\mathcal{O}\left(\left(d_{x}+d_{y}\right) \epsilon^{-4}\right)$}   \\
						C-C  & \cite{liu2020min} & ~ & ~ & ~  \\
						\cline{1-5}
						NC-SC & \multirow{2}{*}{ZO-GDA \cite{wang2020zerothorder}} & \multirow{2}{*}{\tnote{3} $\mathcal{O}\left(\kappa^{5} \epsilon^{-2}\right)$} & \multirow{2}{*}{$\mathcal{O}\left(d_{x}+d_{y}\right)$ }& \multirow{2}{*}{$\mathcal{O}\left(\kappa^{5}\left(d_{x}+d_{y}\right) \epsilon^{-2}\right)$}\\
						UC-C & ~ & ~ & ~ & ~\\
						\cline{1-5}
						NC-SC &  ZO-GDMSA & \multirow{2}{*}{$\mathcal{O}\left(\kappa^{2} \epsilon^{-2} \log \left(\epsilon^{-1}\right)\right)$} & \multirow{2}{*}{$\mathcal{O}\left(d_{x}+ d_{y}\right)$} & \multirow{2}{*}{$\mathcal{O}\left(\kappa\left(d_{x}+\kappa d_{y} \log \left(\epsilon^{-1}\right)\right) \epsilon^{-2}\right)$} \\
						UC-C & \cite{wang2020zerothorder} & ~ & ~ & ~ \\
						\cline{1-5}
						NC-C & \bf{ZO-AGP} &$\mathcal{O}\left(\epsilon^{-4}\right)$ &$\mathcal{O}\left(d_{x}+d_{y}\right)$ &$\mathcal{O}\left(\left(d_{x}+d_{y}\right) \epsilon^{-4}\right)$\\   
						\cline{2-5}              
						C-C & \bf{ZO-BAPG} & $\mathcal{O}\left(\epsilon^{-4}\right)$ & $\mathcal{O}\left(K d_{x}+d_{y}\right)$ & $\mathcal{O}\left(\left(Kd_{x}+d_{y}\right) \epsilon^{-4}\right)$\\
						\hline
					\end{tabular}
					\begin{tablenotes}
						\footnotesize
						\item[1] Total complexity is the product of iteration complexity and per iteration function value  estimation numbers.
						\item[2] ``NC-SC" denotes nonconvex-strongly concave and ``NC-C" denotes nonconvex-concave. \\  ``UC(C)-C" denotes unconstrained (constrained) on $x$ and constrained on $y$.
						\item[3] $\kappa\  denotes \ condition \  number $.
					\end{tablenotes}
			\end{threeparttable}}
			\label{table}
		\end{table}
		
			\begin{remark}\label{remark:2}
			If we set $\rho=\frac{1}{10L_{y}}$,
			From Theorem \ref{thm:main_block},  by plugging $\rho$ and other constants, we can conclude that for any given $\varepsilon \in (0,1)$, $T\leq \max\{(\frac{8sC}{\varepsilon^{2}C_{1}}+2)^{2}, \frac{ 9^{4}\sigma_{y}^{4}L_{y}^{4}}{\varepsilon^{4}}\}$ $= \mathcal{O}\left( \varepsilon^{-4} \right)$. 
			This implies that the total complexity of the proposed ZO-BAPG algorithm to obtain an $\varepsilon$-stationary point for nonconvex-concave minimax problems is bounded by $\mathcal{O}\left(\left(Kd_{x}+d_{y}\right) \epsilon^{-4}\right)$.
		\end{remark}

		We summarize some of the existing zeroth-order algorithms for deterministic nonconvex-(strongly) concave minimax problems in Table \ref{table}.
		
		\section{Numerical Experiments}\label{Experiments}
		
		In this section,  we do some numerical comparisons to solve a data poisoning problem which is a black-box nonconvex minimax problem. All the numerical tests are implemented in Python 3.7 and run in a laptop with 2.70GHz processor, 16GB RAM. 
		
		\subsection{Data poisoning against logistic regression}
		\begin{figure}
			\centering
			\subfigure[]{
				\label{fig_a}
				\includegraphics[width=170pt]{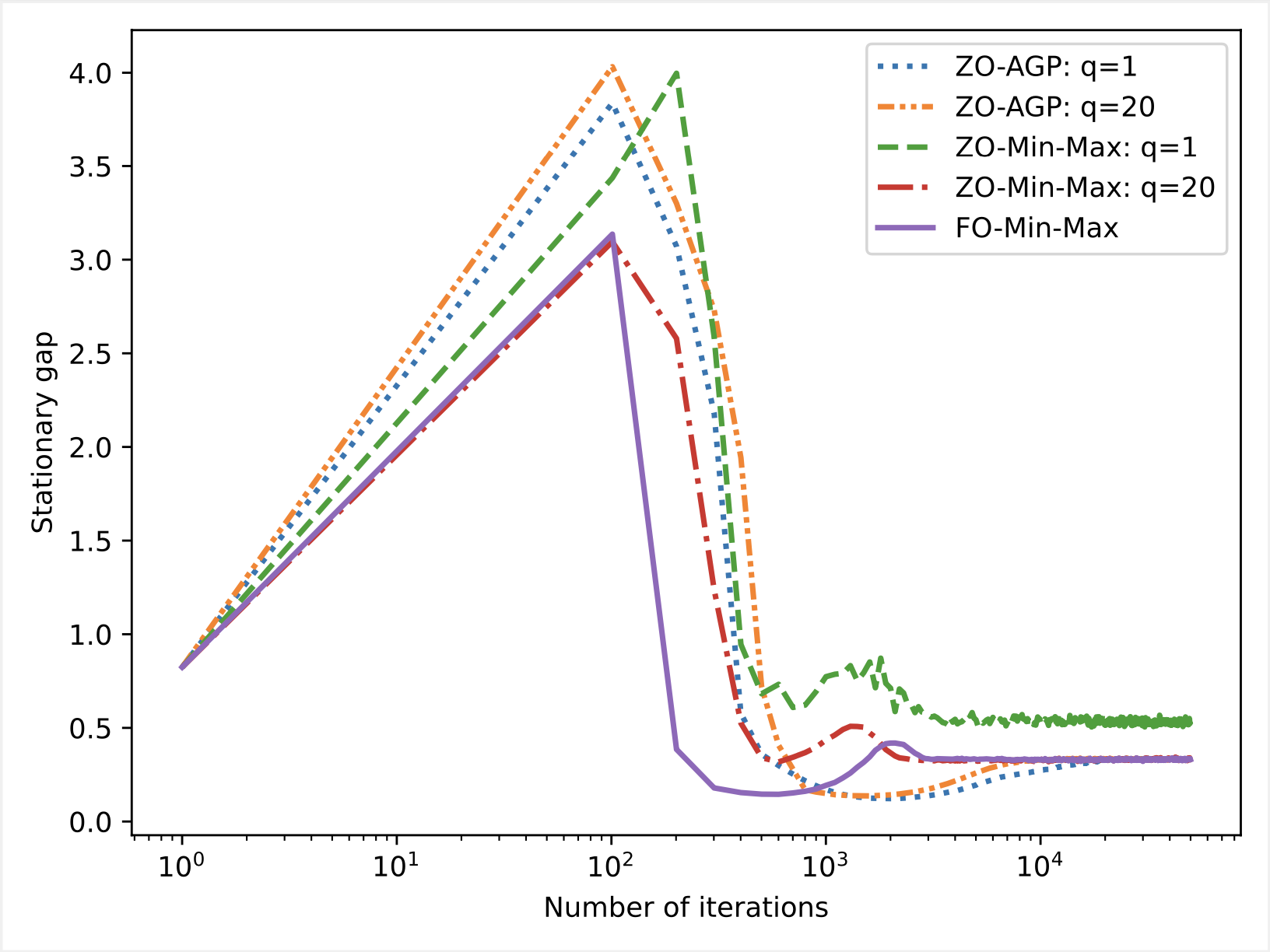}}
			\subfigure[]{
				\label{fig_b}
				\includegraphics[width=170pt]{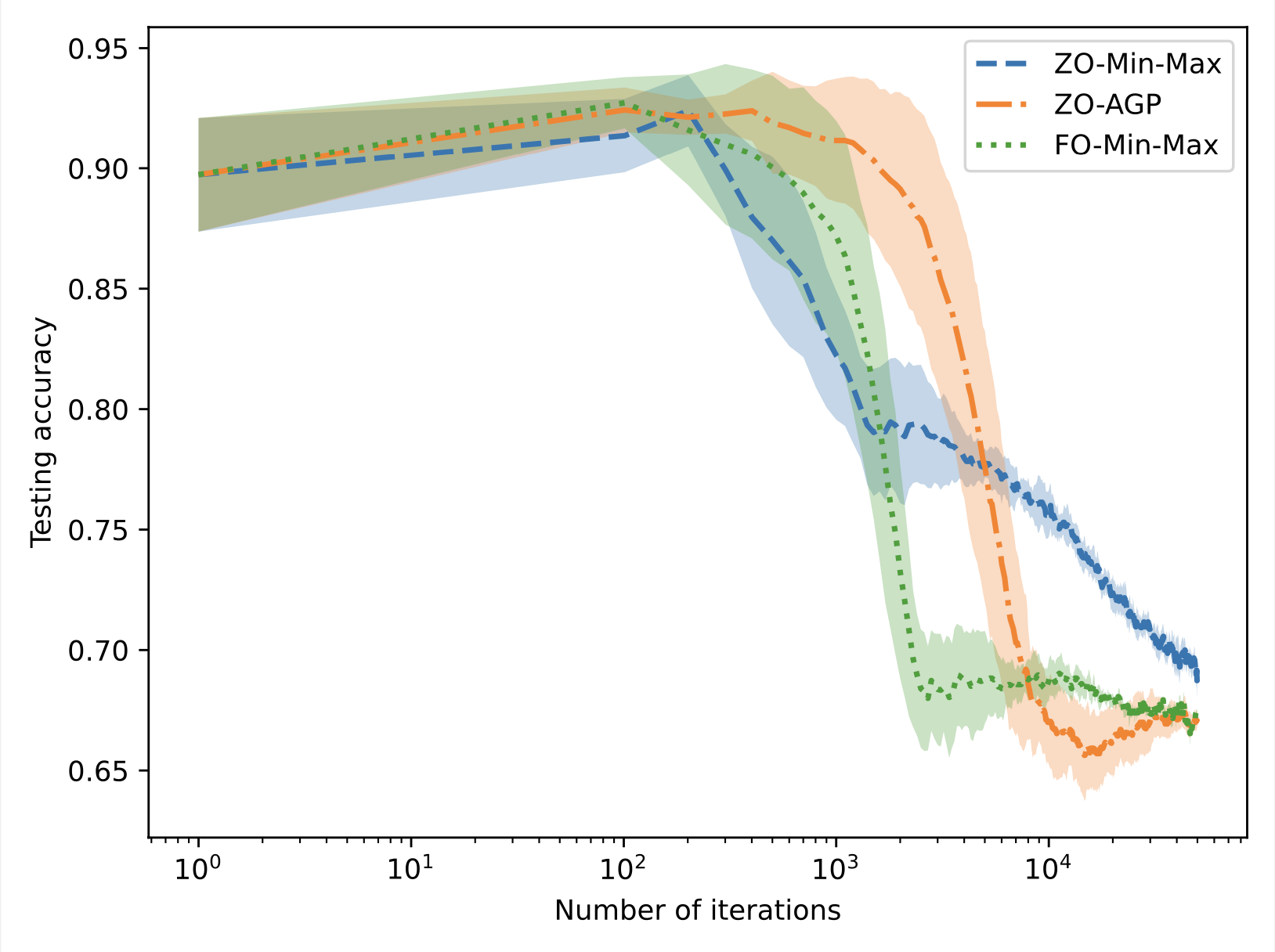}}
			\label{img1}
			\caption{Performance of ZO-Min-Max \cite{liu2020min} and Algorithm \ref{algo1} in data poisoning against logistic regression.}
		\end{figure}
		Data Poisoning \cite{jagielski2018manipulating} is an adversarial attack that the attacker tries to manipulate the training dataset. The goal of the attacker is to corrupt the training dataset so that predictions on the dataset will be modified in the testing phase when using a machine learning model. 
		
		\textbf{Background Setup.} In our experiment, we generate a dataset that contains k=1000 samples  $\left\{z_{i}, t_{i}\right\}_{i=1}^{k}$, where $z_{i} \in \mathbb{R}^{100}$ are sampled from $\mathcal{N}(\mathbf{0}, \mathbf{I})$, and $t_{i}=1$ if $\frac{1}{1+e^{-\left(z_{i}^{T} \boldsymbol{\theta}^{*}+\nu_{i}\right)}} >0.5$, $t_{i}=0$ otherwise, with $\nu_{i} \in \mathcal{N}\left(0,10^{-3}\right)$. Here we choose $ \theta^{*} $ as the base model parameters. The dataset then will be randomly splited into training dataset $ D_{train} $ and testing dataset $ D_{test} $. $\mathcal{D}_{\mathrm{train}}=\mathcal{D}_{\mathrm{tr}, 1} \cup \mathcal{D}_{\mathrm{tr}, 2}, \mathcal{D}_{\mathrm{tr}, 1}$ represents the poisoned dataset which is the subset of the training dataset.
		
		This problem then can be formulated as below:
		\begin{equation}
			\underset{\|x\|_{\infty} \leq \epsilon}{\operatorname{maximize}}\  \underset{\theta}{\operatorname{minimize}} \ f(x, \theta):=F_{\mathrm{tr}}\left(x, \theta ; \mathcal{D}_{train}\right),
			\label{Data_poison}
		\end{equation}
		where $x\in \mathbb{R}^{100}$ and $\theta\in \mathbb{R}^{100}$ are optimization variables. Vector $ x $ is a perturbation which means data in the corrupted dataset is changed to $ z_{i}+x $. $F_{\operatorname{tr}}\left(x, \theta ; \mathcal{D}_{train}\right)=h\left(x, \boldsymbol{\theta} ; \mathcal{D}_{\mathrm{tr}, 1}\right)+h\left(\mathbf{0}, \boldsymbol{\theta} ; D_{\mathrm{tr}, 2}\right)$,  we futher denote
		$$h(x, \boldsymbol{\theta} ; \mathcal{D})=-\frac{1}{\mathcal|{D}|} \sum_{\left(z_{i}, t_{i}\right) \in \mathcal{D}}\left[t_{i} \log \left(h\left(x, \boldsymbol{\theta} ; z_{i}\right)\right)+\left(1-t_{i}\right) \log (1-\right.\left.\left.h\left(x, \boldsymbol{\theta} ; z_{i}\right)\right)\right].$$
		where $h\left(x, \boldsymbol{\theta} ; z_{i}\right)=\frac{1}{1+e^{-\left(z_{i}+x\right)^{T} \boldsymbol{\theta}}}$ which is a logistic regression model for binary classification.  \eqref{Data_poison} can also be written in the form $\min _{x} \max _{\boldsymbol{\theta}}-f(x, \boldsymbol{\theta})$. 
		
		\textbf{Experiment Setup.}
		We compare our algorithm with ZO-Min-Max and FO-Min-Max, which is the first-order formulation of ZO-Min-Max \cite{liu2020min} . Note that \eqref{Data_poison} is different from the experiment in \cite{liu2020min}. \eqref{Data_poison} adopt the ordinary least square setting in \cite{jagielski2018manipulating} with changed loss function. 
		Here we compare all algorithms to solve \eqref{Data_poison} directly, whereas in \cite{liu2020min}, they add a regularization term with $\lambda>0$ as a regularization parameter which makes the problem to be strongly-concave in $ \theta $.
		
		We set the poisoning ratio $ \frac{\left|\mathcal{D}_{\mathrm{tr}, 1}\right|}{\left|\mathcal{D}_{\mathrm{train}}\right|}=10 \%,\epsilon=2$ and iteration $T= 50000$ for the test. Moreover, we set the learning rate $ \alpha_{t} $, $ \beta $,$ \lambda_{t} $ to be $ \frac{5}{100+\sqrt{t}} $, $ 0.02 $, $ \frac{0.1}{t^{1/4}} $ respectively. We set $ \alpha=0.02, \beta=0.05 $ in ZO-Min-Max and FO-Min-Max which is chosen the same as in \cite{liu2020min}. We run 10 independent trials and compare the average performance of three algorithms.
		
		\textbf{Results.}  Figure \ref{fig_a} shows the average stationary gap of the three tested algorithms with different choices of $q$ and the shaded part around the dotted line denotes the standard deviation over 10 independent runs. The stationary gap of ZO-AGP algorithm reaches that of FO-Min-Max algorithm, and both of them perform better than ZO-Min-Max. Figure \ref{fig_b} shows the testing accuracy under the setting of $q=20$. We could see that ZO-AGP yields better accuracy than that of ZO-Min-Max.

		\subsection{Distributed Nonconvex Sparse Principal Component Analysis}
		In this subsection, we compare the proposed ZO-BAPG algorithm with the ZO-Min-Max algorithm for solving a distributed version of nonconvex $\ell_1 $ penalized, nonnegative, sparse principal component analysis (SPCA) problem\cite{Asteris}. 

		\textbf{Background Setup.} In our experiment,
		we construct a connected graph $\mathcal{G}=\{\mathcal{V}, \mathcal{E} \}$ with $\lvert \mathcal{V} \rvert=N$ vertices and $\lvert \mathcal{E} \rvert = E$ edges. Denote $r=\lfloor N/3 \rfloor $.  Define $B \in \mathbb{R}^{E \times N}$ as the incidence matrix, i.e., assuming that the edge $e$ is incident on vertices $i$ and $j$, with $i>j$ we have that if $v=i$, $B_{e v} = 1 $; and if $v=j$, $B_{e v} = -1 $; otherwise $B_{ev}=0$. Let $\otimes$ be the Kronecker product, and $I_d$ is a $d \times d$ identity matrix. According to~\cite{hajinezhad2019perturbed}, the problem can be formulated as
		\begin{align}\label{example4.2}
			\min_{x} \max_{y} & - \sum_{i=1}^{N} (x^{(i)})^\top \Sigma_i x^{(i)} + \sum_{i=1}^{r} \frac{N\mu}{r}\Vert x^{(i)} \Vert_1 + y^\top (B \otimes I_d)x \nonumber \\
			s.t. \quad & \|x^{(i)} \|^2 \leq 1, \quad i=r+1,\cdots 2r \\
			& x^{(i)} \geq 0, \quad i=2r+1,\cdots N.\nonumber
		\end{align}
		where $x=[x^{(1)}; \cdots; x^{(N)}] \in \mathbb{R}^{Nd \times 1}$ with $x^{(i)}\in \mathbb{R}^{d \times 1} $,  $\Sigma_i \in \mathbb{R}^{d \times d}$ is the covariance matrix for the mini-batch data in node $i$, 
		$\mu > 0$ is a constant that controls the sparsity.
		
		\textbf{Experiment Setup.}
		The graph $\mathcal{G}$ is generated based on the scheme proposed in~\cite{Yildiz} that $N$ nodes uniformly distributed over a unit square and two nodes connect to each other if their distance is less than $R$. Set  $N=10$, $R=0.7$, $d=8$, and $\mu=0.01$.  The parameters are chosen as $ \beta=0.05$, $ \lambda_{t}=\frac{0.01}{t^{1/4}} $, $\gamma_t=1000$, $\tau_t=100t^{1/2}$ in ZO-BAPG and $ \alpha=0.01, \beta=0.05 $ in ZO-Min-Max. Let $\|(B \otimes I_d)x\|^2$ denote the constraint violation (Cons-Vio). The smaller the Cons-Vio, the better the feasibility of the obtained solution.
		
\begin{figure}
			\centering
			\subfigure[]{
				\label{fig_gap}
				\includegraphics[width=170pt]{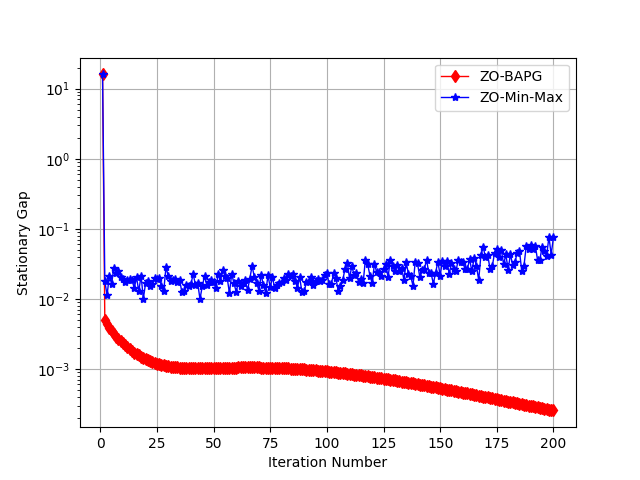}}
			\subfigure[]{
				\label{fig_cv}
				\includegraphics[width=170pt]{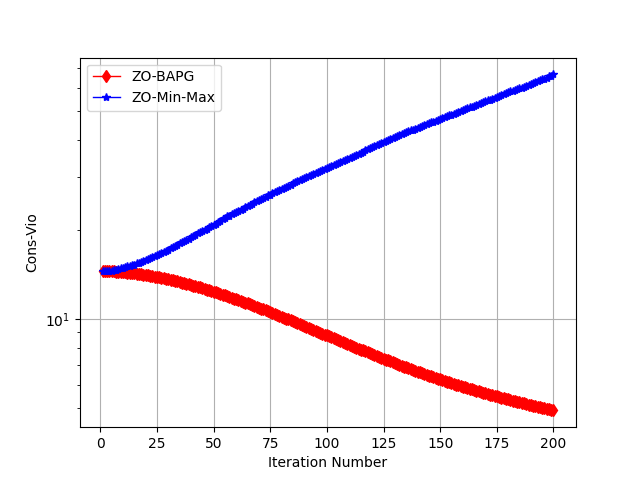}}
			\caption{Performance of ZO-Min-Max \cite{liu2020min} and Algorithm \ref{algo2} in distributed nonconvex quadratic problem.}
		\end{figure}

			\textbf{Results.} Figure \ref{fig_gap} and  Figure \ref{fig_cv} show the stationary gap and the constraint violation degree of the ZO-BAPG algorithm and the ZO-Min-Max algorithm for solving \eqref{example4.2}. It can be seen from Figure \ref{fig_gap} and \ref{fig_cv} that as the number of iterations increases, neither the stationary gap nor the constraint violation degree of the iteration points obtained by the ZO-Min-Max algorithm decreases, which means that the limit of the iterative point sequence obtained by the ZO-Min-Max algorithm is not a feasible solution of the original SPCA problem. While both the stationary gap and the constraint violation degree of the ZO-BAPG algorithm decrease as the number of iterations increases, which shows that the proposed ZO-BAPG algorithm has better numerical performance than the ZO-Min-Max algorithm in solving problem \eqref{example4.2}.

		\section{Conclusion}\label{section_conclu}
		This paper proposes two zeroth-order algorithms for solving nonconvex-concave minimax problems and their comlexity are also analyzed. To the best of our knowledge, this is the first time that zeroth-order algorithms with iteration complexity gurantee are developed for solving both general smooth and block-wise nonsmooth nonconvex-concave minimax problems. Furthermore, the effectiveness and efficiency of the proposed ZO-AGP algorithm is verified by some limited numerical tests on data poisoning attack problem  and  distributed nonconvex sparse principal component analysis problem. 
		
		%

		\section*{Acknowledgments}
We would like to thank two anonymous reviewers for their helpful comments and suggestions to improve this paper.
		
		\bibliographystyle{siamplain}
		
	\end{document}